\documentclass[fleqn]{article}
\usepackage{amssymb,latexsym,theorem}

\newcommand{\C}{\mathbb{C}}

\newcommand{\PG}{\mathrm{PG}}

\newcommand{\F}{\mathbb{F}}

\newtheorem{lemma}{Lemma}[section]
\newtheorem{prop}[lemma]{Proposition}
\newtheorem{theo}[lemma]{Theorem}
\newtheorem{co}[lemma]{Corollary}

\newtheorem{theorem}{Theorem}

\newtheorem{Lemma}[theorem]{Lemma}

\def\pr{\noindent{\bf Proof. }}
\def\eop{\hspace*{\fill}$\Box$}

\title{Veronesean embeddings of dual polar spaces of orthogonal type}
\author{Ilaria Cardinali and Antonio Pasini}
\date{}

\begin{document}
\maketitle
\begin{abstract}
Given a point-line geometry $\Gamma$ and a pappian projective space $\cal S$, a veronesean embedding of $\Gamma$ in $\cal S$ is an injective map $e$ from the point-set of $\Gamma$ to the set of points of $\cal S$ mapping the lines of $\Gamma$ onto non-singular conics of $\cal S$ and such that $e(\Gamma)$ spans $\cal S$. In this paper we study veronesean embeddings of the dual polar space $\Delta_n$ associated to a non-singular quadratic form $q$ of Witt index $n \geq 2$ in $V = V(2n+1,\mathbb{F})$. Three such embeddings are considered, namely the Grassmann embedding $\varepsilon^{\mathrm{gr}}_n$ which maps a maximal singular subspace $\langle v_1,..., v_n\rangle$ of $V$ (namely a point of $\Delta_n$) to the point $\langle \wedge_{i=1}^nv_i\rangle$ of $\mathrm{PG}(\bigwedge^nV)$, the composition $\varepsilon^{\mathrm{vs}}_n := \nu_{2^n}\circ \varepsilon^{\mathrm{spin}}_n$ of the spin (projective) embedding $\varepsilon^{\mathrm{spin}}_n$ of $\Delta_n$ in $\mathrm{PG}(2^n-1,\mathbb{F})$ with the quadric veronesean map $\nu_{2^n}:V(2^n,\mathbb{F})\rightarrow V({{2^n+1}\choose 2}, \mathbb{F})$, and a third embedding $\tilde{\varepsilon}_n$ defined algebraically in the Weyl module $V(2\lambda_n)$, where $\lambda_n$ is the fundamental dominant weight associated to the $n$-th simple root of the root system of type $B_n$. We shall prove that $\tilde{\varepsilon}_n$ and $\varepsilon^{\mathrm{vs}}_n$ are isomorphic. If $\mathrm{char}(\F)\neq 2$ then $V(2\lambda_n)$ is irreducible and $\tilde{\varepsilon}_n$ is isomorphic to $\varepsilon^{\mathrm{gr}}_n$ while if $\mathrm{char}(\F)= 2$ then $\varepsilon^{\mathrm{gr}}_n$ is a proper quotient of $\tilde{\varepsilon}_n$. In this paper we shall study some of these submodules. Finally we turn to universality, focusing on the case of $n = 2$. We prove that if $\F$ is a finite field of odd order $q > 3$ then $\varepsilon^{\mathrm{sv}}_2$ is relatively universal. On the contrary, if $\mathrm{char}(\F)= 2$ then $\varepsilon^{\mathrm{vs}}_2$ is not universal. We also prove that if $\F$ is a perfect field of characteristic 2 then $\varepsilon^{\mathrm{vs}}_n$ is not universal, for any $n \geq 2$.
\end{abstract}

{\scriptsize
}

\section{Introduction}\label{Introduction}

\subsection{Projective and veronesean embeddings}\label{Intro:Emb}

We firstly recall a few basics on projective embeddings. A {\em projective embedding} of a point-line geometry $\Gamma = ({\cal P},{\cal L})$ in the projective space $\PG(V)$ of a vector space $V$ is an injective mapping $\varepsilon$ from the point-set $\cal P$ of $\Gamma$ to the set of points of $\PG(V)$ such that $\varepsilon$ maps every line of $\Gamma$ surjectively onto a line of $\PG(V)$ and $\varepsilon({\cal P})$ spans $\PG(V)$. Henceforth we will freely commit the abuse of regarding $V$ as the codomain of $\varepsilon$ instead of $\PG(V)$, thus writing $\varepsilon:\Gamma\rightarrow V$ instead of $\varepsilon:\Gamma\rightarrow\PG(V)$. Accordingly, if $p\in {\cal P},$ we regard $\varepsilon (p)$ as a $1$-dimensional subspace of $V$ and we take the dimension of $V$ as the {\em dimension} $\mathrm{dim}(\varepsilon)$ of $\varepsilon$.

If $\F$ is the underlying division ring of $V$ then we say that $\varepsilon$ is {\em defined over} $\F$, also that $\varepsilon$ is a {\em projective $\F$-embedding} for short.

Given two projective $\F$-embeddings $\varepsilon_1:\Gamma\rightarrow V_1$ and $\varepsilon_2:\Gamma\rightarrow V_2$, a {\em morphism} $f:\varepsilon_1\rightarrow \varepsilon_2$ from $\varepsilon_1$ to $\varepsilon_2$ is a semi-linear mapping $f:V_1\rightarrow V_2$ such that $\varepsilon_2 = f\cdot\varepsilon_1$. Note that, since $\langle \varepsilon_2({\cal P})\rangle = V_2$, the equality $\varepsilon_2 = f\cdot \varepsilon_1$ forces $f:V_1\rightarrow V_2$ to be surjective. If $f$ is bijective then $f$ is called an {\em isomorphism}. When $\varepsilon_1$ and $\varepsilon_2$ are isomorphic we write $\varepsilon_1\cong \varepsilon_2$. Note that if a morphism $f:\varepsilon_1\rightarrow\varepsilon_2$ exists then $f$ is uniquely determined modulo scalars. If a morphism exists from $\varepsilon_1$ to $\varepsilon_2$ then we write $\varepsilon_1\geq \varepsilon_2.$  When $\varepsilon_1 \geq \varepsilon_2$ but $\varepsilon_1\not\cong \varepsilon_2$ we write $\varepsilon_1 > \varepsilon_2$.

If $f:V_1\rightarrow V_2$ is a morphism from $\varepsilon_1$ to $\varepsilon_2$ then $V_2 \cong V_1/\mathrm{ker}(f)$ and an embedding $\varepsilon_1/\mathrm{ker}(f):\Gamma\rightarrow \PG(V_1/\mathrm{ker}(f))$ can be defined mapping every point $p\in{\cal P}$ of $\Gamma$ to $\langle \varepsilon_1(p), \mathrm{ker}(f)\rangle/\mathrm{ker}(f)$. We say that $\varepsilon_1/\mathrm{ker}(f)$ is a {\em quotient} of $\varepsilon_1$. Clearly, if $f:V_1\rightarrow V_2$ is a morphism from $\varepsilon_1$ to $\varepsilon_2$ then $\varepsilon_2\cong \varepsilon_1/\mathrm{ker}(f)$. In view of this fact, we take the liberty to call $\varepsilon_2$ a {\em quotient} of $\varepsilon_1$ (a {\em proper quotient} if $\varepsilon_1\not\cong\varepsilon_2$).  We also call the morphism $f:\varepsilon_1\rightarrow\varepsilon_2$ the {\em projection} of $\varepsilon_1$ onto $\varepsilon_2$.

Following Kasikova and Shult \cite{KS}, we say that a projective embedding of a point-line geometry $\Gamma$ is {\em relatively universal} when it is not a proper quotient of any other
projective embedding of $\Gamma$. Every projective embedding $\varepsilon$ of $\Gamma$ admits a {\em hull} $\tilde{\varepsilon}$, uniquely determined up to isomorphism by the following property: $\tilde{\varepsilon}$ is a projective embedding of $\Gamma$, $\varepsilon$ is a quotient of $\tilde{\varepsilon}$ and we have $\tilde{\varepsilon}\geq \varepsilon'$ for every projective embedding $\varepsilon'$ of $\Gamma$ such that $\varepsilon'\geq \varepsilon$ (see Ronan \cite{Ronan}, where a construction of the hull of a projective embedding by means of
a suitable presheaf of 1- and 2-dimensional $\F$-vector spaces is also given). Clearly, the hull $\tilde{\varepsilon}$ of $\varepsilon$ is relatively universal. A projective embedding is relatively universal if and only if it is its own hull.

A projective $\F$-embedding $\tilde{\varepsilon}$ of $\Gamma$ is {\em absolutely universal} if all projective $\F$-embeddings of $\Gamma$ are quotients of $\tilde{\varepsilon}$. The absolutely universal projective $\F$-embedding of $\Gamma$, if it exists, is uniquely determined up to isomorphisms. It is the hull of all projective $\F$-embeddings of $\Gamma$. Obviously, every absolutely universal projective embedding is relatively universal. If $\Gamma$ admits the absolutely universal projective $\F$-embedding then the converse also holds true: all relatively universal projective $\F$-embeddings of $\Gamma$ are absolutely universal.

Given an embedding $\varepsilon:\Gamma\rightarrow V$ and an automorphism $g$ of $\Gamma$, a {\em lifting} of $g$ through $\varepsilon$ is a semi-linear mapping $\varepsilon(g):V\rightarrow V$ such that $\varepsilon(g)\cdot \varepsilon = \varepsilon\cdot g$. The lifting $\varepsilon(g)$ of $g$, if it exists, is uniquely determined modulo scalars. Clearly, it is invertible. Given a group $G$ acting on $\Gamma$ as a group of automorphisms, the embedding $\varepsilon$ is said to be $G$-{\em homogeneous} if for every $g\in G$ the automorphism of $\Gamma$ induced by $g$ lifts through $\varepsilon$ to a semi-linear map of $V$. It is easy to see that if $\varepsilon$ is relatively universal then it is $\mathrm{Aut}(\Gamma)$-homogeneous.

We now turn to veronesean embeddings. Veronesean embeddings are defined just like projective embeddings except that they map lines onto conics. Explicitly, a {\em veronesean embedding} of a point-line geometry $\Gamma = ({\cal P},{\cal L})$ in (the projective space $\PG(V)$ of) a vector space $V$ defined over a commutative division ring (namely a field) is an injective mapping $\varepsilon$ from the point-set $\cal P$ of $\Gamma$ to the set of points of $\PG(V)$ such that $\varepsilon$ maps every line of $\Gamma$ onto a non-singular conic of $\PG(V)$, $\varepsilon(p)\cap \langle \varepsilon(l)\rangle =0$ for every non-incident point-line pair $(p,l)$ and $\varepsilon({\cal P})$ spans $\PG(V)$. All definitions and conventions stated for projective embeddings can be rephrased for veronesean embeddings word by word. A few minor modifications are needed only in the construction of the hull of a veronesean embedding by means of a presheaf: now we need a presheaf of 1- and 3-dimensional vector spaces (see \cite{Pas}, where the hulls we are interested in here are called {\em linear hulls}).

When $\F$ is a perfect field of characteristic 2 we can also define morphisms from a veronesean $\F$-embedding $\varepsilon_1:\Gamma\rightarrow V_1$ to a projective $\F$-embedding $\varepsilon_2:\Gamma\rightarrow V_2$. Such a morphism is a semi-linear mapping $f:V_1\rightarrow V_2$ such that for every line $l$ of $\Gamma$, $\mathrm{ker}(f)\cap\langle\varepsilon_1(l)\rangle$ is the nucleus $n_l$ of the conic $\varepsilon_1(l)$ and $\varepsilon_2(l)$ is canonically isomorphic to the set of lines of $\langle \varepsilon_1(l)\rangle$ through $n_l$. Note that each of those lines meets $\varepsilon_1(l)$ in a point, since $\F$ is assumed to be perfect.

A situation like the above can also be considered when $\F$ is a non-perfect field of characteristic 2, but in that case $\varepsilon_2$ is a {\em lax} embedding, where the image $\varepsilon_2(l)$ of a line $l$ of $\Gamma$ is a possibly proper subset of a line of $\PG(V_2)$. However, we prefer to keep lax embeddings out of the scope of this paper.

\subsection{The geometries and the embeddings to be considered in this paper}\label{1.2}

Let $V := V(2n+1,\F)$ be a $(2n+1)$-dimensional vector space over a field $\F$ and $q$ a non singular quadratic form of $V$ of Witt index $n\geq 2$. Let $\Delta$ the building of type $B_n$ where the elements of type $k = 1, 2,..., n$ ($k$-elements for short) are the $k$-dimensional subspaces of $V$ totally singular for $q$, with containment as the incidence relation.

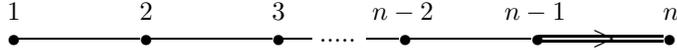
\begin{picture}(310,36)(0,0)
\put(20,8){$\bullet$}
\put(23,11){\line(1,0){47}}
\put(70,8){$\bullet$}
\put(73,11){\line(1,0){47}}
\put(120,8){$\bullet$}
\put(123,11){\line(1,0){12}}
\put(138,10){$.....$}
\put(156,11){\line(1,0){12}}
\put(168,8){$\bullet$}
\put(171,11){\line(1,0){47}}
\put(218,8){$\bullet$}
\put(241,8){\large{$>$}}
\put(221,10){\line(1,0){47}}
\put(221,12){\line(1,0){47}}
\put(268,8){$\bullet$}
\put(20,18){1}
\put(70,18){2}
\put(120,18){3}
\put(158,18){$n-2$}
\put(208,18){$n-1$}
\put(268,18){$n$}
\end{picture}

\noindent
For $1\leq k\leq n$, the $k$-shadow of a flag $F$ of $\Delta$ is the set of $k$-elements incident to $F$. The $k$-{\em grassmannian} $\Delta_k$ of $\Delta$ is the point-line geometry defined as follows. The points of $\Delta_k$ are the $k$-elements of $\Delta$. When $1< k < n$ the lines of $\Delta_k$ are the $k$-shadows of the flags of $\Delta$ of type $\{k-1,k+1\}$. The lines of $\Delta_1$ are the 1-shadows of the 2-elements of $\Delta$. The lines of $\Delta_n$ are the $n$-shadows of the $(n-1)$-elements of $\Delta$. For such an element $X$, let $l_X$ be its $n$-shadow. Then
\begin{equation}\label{lines}
l_{X} = \{Z\mid X \subset Z\subset X^{\perp} , ~\mathrm{dim}(Z) = n, ~ Z ~\mbox{totally singular}\}
\end{equation}
where $X^\perp$ is the orthogonal of $X$ with respect to $q$. (Recall that $X$ is an $(n-1)$-dimensional totally singular subspace of $V$.) The vector space $X^\perp/X$ is 3-dimensional and $l_{X}$ is a non-singular conic in the projective plane $\PG(X^\perp/X)$. The geometry $\Delta_n$ is called a {\em dual polar space} of type $B_n$, while $\Delta_1$ is the polar space associated to the building $\Delta$.

In this paper we are mainly interested in the dual polar space $\Delta_n$ and its veronesean embeddings. Let $W_n :=\bigwedge^n V$. We recall that $\mathrm{dim}(W_n) = {{2n+1}\choose n}$. The points of $\Delta_n$ are the $n$-dimensional totally singular subspaces of $V$. The {\em Grassmann embedding} $e^{\mathrm{gr}}_n$ of $\Delta_n$ maps every such subspace $\langle v_1,v_2,\dots, v_n\rangle$ to the point $\langle v_1\wedge v_2\wedge\dots \wedge v_n \rangle$ of $\PG(W_n)$. If $X$ is an $(n-1)$-element of $\Delta$ then the set of points $\varepsilon_n^{\mathrm{gr}}(l_X) = \{\varepsilon_n^{\mathrm{gr}}(Y)\}_{Y\in l_X}$ is a non-degenerate conic of $\PG(W_n)$. So, $\varepsilon_n^{\mathrm{gr}}$ is a veronesean embedding of $\Delta_n$ in the subspace $\langle\varepsilon_n^{\mathrm{gr}}(\Delta_n)\rangle$ of $W_n$ spanned by $\varepsilon_n^{\mathrm{gr}}(\Delta_n)$, where we take the liberty of using the symbol $\Delta_n$ to denote also the point-set of the geometry $\Delta_n$. As proved in \cite{CP1}, if $\mathrm{char}(\F)\neq 2$ then $\langle\varepsilon_n^{\mathrm{gr}}(\Delta_n)\rangle = W_n$ while if $\mathrm{char}(\F) = 2$ then $\langle\varepsilon_n^{\mathrm{gr}}(\Delta_n)\rangle$ is a subspace of $W_n$ of codimension equal to ${{2n+1}\choose{n-2}}$.

The dual polar space $\Delta_n$ also admits a projective embedding, namely the {\em spin embedding} $\varepsilon^{\mathrm{spin}}_n : \Delta_n \rightarrow V(2^n,\F)$. We refer the reader to Buekenhout and Cameron \cite{BC} for a concise description of this embedding. It is worth mentioning that when $\mathrm{char}(\F) \neq 2$ the embedding $\varepsilon_n^{\mathrm{spin}}$ is relatively universal (Blok and Brouwer \cite{BB}; also Cooperstein and Shult \cite{CS}). Hence it is absolutely universal, since $\Delta_n$ admits the absolutely universal embedding (Kasikova and Shult \cite{KS}).

Let $\nu_{2^n}$ be the usual quadric veronesean map from $V(2^n,{\mathbb{F}})$ to $V({{2^n+1}\choose 2},{\mathbb{F}})$, which maps $(x_1,...,x_{2^n})\in V(2^n,{\mathbb{F}})$ to the vector
\[(x_1^2,\dots, x_{2^n}^2,x_1x_2,\dots, x_1x_{2^n},x_2x_3,\dots, x_2x_{2^n}, \dots, x_{2^n-1}x_{2^n}).\]
The mapping $\nu_{2^n}$ defines a veronesean embedding of the point-line geometry $\PG(2^n-1,\F)$ in $V({{2^n+1}\choose 2},\F)$, which we also denote by the symbol $\nu_{2^n}$. The composition $\varepsilon^{\mathrm{vs}}_n := \nu_{2^n}\cdot\varepsilon^{\mathrm{spin}}_n$ is a veronesean embedding of $\Delta_n$ in a subspace of $V({{2^n+1}\choose 2},\F)$. We call it the {\em veronesean-spin embedding} of $\Delta_n$.

Before describing the third veronesean embedding of $\Delta_n$ we must fix some notation for groups. Throughout this paper $G := \mathrm{Spin}(2n+1,\F)$ is the spin group of rank $n$ defined over $\F$. Thus, $G$ is the universal Chevalley group of type $B_n$ defined over $\F$. The adjoint group of type $B_n$ is $\overline{G} := \mathrm{SO}(2n+1,\F)$ ($= \mathrm{PSO}(2n+1,\F)$). We recall that if $\mathrm{char}(\F)\neq 2$ then $G$ is a non-split central extension of $\overline{G}$ by a group of order 2 while if $\mathrm{char}(\F) = 2$ then $G = \overline{G}$.

Each of the embeddings $\varepsilon^{\mathrm{gr}}_n$, $\varepsilon^{\mathrm{spin}}_n$ and $\varepsilon_n^{\mathrm{vs}}$ considered so far is $G$-homogeneous. The vector space $V(2^n,\F)$, regarded as a $G$-module via $\varepsilon^{\mathrm{spin}}_n$, is called the {\em spin module}. We shall denote it by the symbol $W^{\mathrm{spin}}_n.$ We call the codomains $\langle \varepsilon_n^{\mathrm{gr}}(\Delta_n)\rangle$ and $\langle \varepsilon_n^{\mathrm{vs}}(\Delta_n)\rangle$ of $\varepsilon_n^{\mathrm{gr}}$ and $\varepsilon_n^{\mathrm{vs}}$ the {\em grassmann module} and the {\em veronese-spin module} respectively and we denote them by the symbols $W_n^{\mathrm{gr}}$ and $W_n^{\mathrm{vs}}$. As remarked before, $W^{\mathrm{gr}}_n = W_n$ when $\mathrm{char}(\F) \neq 2$ while if $\mathrm{char}(\F) = 2$ then $W_n^{\mathrm{gr}}$ is a proper submodule of $W_n$. When $\mathrm{char}(\F) \neq 2$ the group $G$ acts as $\overline{G}$ on $W^{\mathrm{gr}}_n$ ($= W_n$) as well as in $W^{\mathrm{vs}}_n$, but it acts faithfully in $W^{\mathrm{spin}}_n.$ Thus, $W^{\mathrm{spin}}_n$ is not a $\overline{G}$-module.

We now turn to the third veronesean embedding of $\Delta_n$. Let $\lambda_1, \lambda_2,..., \lambda_n$ be the fundamental dominant weights for the root system of type $B_n$, numbered in the usual way (see the picture at the beginning of this subsection). For $\lambda = \lambda_n$ or $\lambda = 2\cdot\lambda_n$, let $V(\lambda)$ be the Weyl $G$-module with highest weight $\lambda$. An embedding $\varepsilon_\lambda$ of $\Delta_n$ into $V(\lambda)$ can be created as follows. Let $v_0$ be a highest weight vector of $V(\lambda)$. Then the $G$-orbit of $\langle v_0\rangle$ corresponds to the set of points of $\Delta_n$ and, if $P_n$ is the minimal fundamental parabolic subgroup of $G$ of type $n$ and $L_0$ is the $P_n$-orbit of $\langle v_0\rangle$, then the $G$-orbit of $L_0$ corresponds to the set of lines of $\Delta_n$. If $X$ is the point of $\Delta_n$ corresponding to $\langle v_0\rangle$, then  $\varepsilon_{\lambda}$ maps $g(X)$ to $g(\langle v_0\rangle )$, for every $g\in G.$ It is well known that $\varepsilon_{\lambda_n} \cong \varepsilon_n^{\mathrm{spin}}$, namely $V(\lambda_n)\cong W^{\mathrm{spin}}_n.$ On the other hand, $\varepsilon_{2\lambda_n}$ is veronesean, as one can see by computing $L_0$ explicitly. We denote $\varepsilon_{2\lambda_n}$ by the symbol $\tilde{\varepsilon}_n$ and we call it the {\em veronesean Weyl embedding} of $\Delta_n$.

We have $\mathrm{dim}(V(2\lambda_n)) = {{2n+1}\choose n}$, as one can check by using the Weyl dimension formula (see e.g. Humphreys \cite[24.3]{Humphreys}). Hence $\mathrm{dim}(\tilde{\varepsilon}_n) = {{2n+1}\choose n}$.

It is known that $\varepsilon_n^{\mathrm{gr}}\leq \tilde{\varepsilon}_n$. More explicitly, if $\mathrm{char}(\F) \neq 2$ then $\varepsilon_n^{\mathrm{gr}}\cong \tilde{\varepsilon}_n$, namely $V(2\lambda_n) \cong W^{\mathrm{gr}}_n$ ($= W_n$). If $\mathrm{char}(\F) = 2$ then $\varepsilon_n^{\mathrm{gr}}$ is a proper quotient of $\tilde{\varepsilon}_n$ (see \cite{CP1}), namely $W^{\mathrm{gr}}_n$ is a proper quotient of $V(2\lambda_n)$. In this case, if $\pi:V(2\lambda_n)\rightarrow W^{\mathrm{gr}}_n$ is the projection of $\tilde{\varepsilon}_n$ onto $\varepsilon_n^{\mathrm{gr}}$ then $\mathrm{dim}(\mathrm{ker}(\pi)) = {{2n+1}\choose{n-2}}$.

Sometimes in this paper we must consider also the {\em Grassmann embedding} $\varepsilon_k^{\mathrm{gr}}$ and the {\em Weyl embedding} $\tilde{\varepsilon}_k$ of $\Delta_k$ for some $k < n$. They are defined in the same way as $\varepsilon_n^{\mathrm{gr}}$ and $\tilde{\varepsilon}_n$ but for replacing $W_n$ with $W_k := \bigwedge^kV$ and $V(2\lambda_n)$ with $V(\lambda_k)$. The embeddings $\varepsilon_k^{\mathrm{gr}}$ and $\tilde{\varepsilon}_k$ are projective and $\mathrm{dim}(\tilde{\varepsilon}_k) = {{2n+1}\choose k}$. Clearly, $\varepsilon_1^{\mathrm{gr}} = \tilde{\varepsilon}_1$. The embedding $\varepsilon_1^{\mathrm{gr}}$ is absolutely universal (Tits \cite[chapter 8]{Tits}). Let $1 < k < n$. In this case, if $\mathrm{char}(\F) \neq 2$ then $\varepsilon_k^{\mathrm{gr}} = \tilde{\varepsilon}_k$ while if $\mathrm{char}(\F) = 2$ then $\varepsilon_k^{\mathrm{gr}}$ is a proper quotient of $\tilde{\varepsilon}_k$ (see \cite{CP1}). It is worth mentioning that when $\F$ is either a perfect field of positive characteristic or a number field and either $k = 2 < n$ or $k = 3 < n$ then $\tilde{\varepsilon}_k$ is absolutely universal \cite{CP1}.

When $\F$ is a perfect field of characteristic 2 the building $\Delta$ is isomorphic to the building of type $C_n$ associated to a non-degenerate alternating form $\alpha$ on $\overline{V} := V(2n,\F)$, the elements of $\Delta$ of type $k$ being now regarded as $k$-subspaces of $\overline{V}$ totally isotropic for $\alpha$. Thus, we can also define a projective embedding $\varepsilon_k^{\mathrm{sp}}$ of $\Delta_k$ in a subspace $W^{\mathrm{sp}}_k:=\langle \varepsilon^{\mathrm{sp}}_k(\Delta_k)\rangle$ of $\overline{W}_k := \bigwedge^k\overline{V}$,  which maps every totally isotropic $k$-subspace $\langle v_1,..., v_k\rangle$ of $\overline{V}$ onto the point $\langle v_1\wedge...\wedge v_k\rangle$ of $\PG(\overline{W}_k)$. We warn that the embedding $\varepsilon_k^{\mathrm{sp}}$ is projective for $k = n$ too. If $k > 1$ then $W^{\mathrm{sp}}_k$ is a proper subspace of $\overline{W}_k$. In fact $\mathrm{dim}(\varepsilon^{\mathrm{sp}}_k) = {{2n}\choose k}-{{2n}\choose{k-2}}$ while $\mathrm{dim}(\overline{W}_k) = {{2n}\choose k}$.

Let $k = n$. The spin embedding $\varepsilon^{\mathrm{spin}}_n$ is a quotient of $\varepsilon_n^{\mathrm{sp}}$ (Blok, Cardinali and De Bruyn~\cite{BCDB09}, see also Cardinali and Lunardon~\cite{CL}). Actually  $\varepsilon^{\mathrm{spin}}_n < \varepsilon^{\mathrm{sp}}_n$, since
$\mathrm{dim}(\varepsilon_n^{\mathrm{sp}}) = {{2n}\choose n}-{{2n}\choose{n-2}} > 2^n = \mathrm{dim}(\varepsilon^{\mathrm{spin}}_n)$. When $2 < |\F| < \infty$ we can say more: in that case $\varepsilon^{\mathrm{sp}}_n$ is the hull of $\varepsilon_n^{\mathrm{spin}}$ (Cooperstein \cite{Co:1}). Finally, let $n = 2$. In this case $\varepsilon^{\mathrm{sp}}_1 \cong \varepsilon_2^{\mathrm{spin}} < \varepsilon_2^{\mathrm{sp}}\cong \varepsilon_1^{\mathrm{gr}}\cong \tilde{\varepsilon}_1$.

\subsection{Results}\label{sec1.3}

We have previously remarked that $\varepsilon_n^{\mathrm{gr}}\leq\tilde{\varepsilon}_n$ and $\varepsilon_n^{\mathrm{gr}}\cong \tilde{\varepsilon}_n$ if and only if $\mathrm{char}(\F) \neq 2$. The following will be proved in Section \ref{sec2}.

\begin{theorem}\label{main theorem 1}
We have $\varepsilon_n^{\mathrm{vs}}\cong \tilde{\varepsilon}_n$ for any $n\geq 2$ and any choice of $\F$.
\end{theorem}
In other words, $V(2\lambda_n)$ and $W^{\mathrm{vs}}_n$ are isomorphic as $G$-modules. The $G$-module $V(2\lambda_n)$ ($\cong W^{\mathrm{vs}}_n$) is irreducible if and only if $\mathrm{char}(\F)\neq 2$ (see \cite{CP1}). In the sequel we shall focus on the case of $\mathrm{char}(\F) = 2$, but before turning to that case we need to recall a few properties of quadric veronesean maps.

For an integer $d \geq 2$ let $\nu_{d+1}$ be the veronesean embedding of $\PG(d,\F)$ in $V({{d+2}\choose 2},\F)$ induced by the quadric veronesean map from $V(d+1,\F)$ to $V({{d+2}\choose 2},\F)$. Clearly, $\nu_{d+1}$ is $\mathrm{GL}(d+1,\F)$-homogeneous.

Let now $\mathrm{char}(\F) = 2$. The {\em nucleus subspace} of $V({{d+2}\choose 2},\F)$ relative to $\nu_{d+1}$ is the subspace $\cal N$ of $V({{d+2}\choose 2},\F)$ spanned by the nuclei of the conics $\nu_{d+1}(l)$, for $l$ a line of $\PG(d,\F)$ (Thas and Van Maldeghem \cite{TVM2004}). We have ${\cal N}\cap \langle \nu_{d+1}(l)\rangle=n_l$ for every line $l$ of $\PG(d,\F).$ Hence ${\cal N}\cap \nu_{d+1}(p)=0$ for every point $p$ of $\PG(d,\F),$ namely $\mathrm{dim}(\langle {\cal N}\cup \nu_{d+1}(p)\rangle / {\cal N})=1.$ The nucleus subspace $\cal N$ is stabilized by $\mathrm{GL}(d+1,\F)$ in its action on $V({{d+2}\choose 2},\F)$.

Put $d = 2^n-1$. Let $\cal N$ be the nucleus subspace of $V({{d+2}\choose 2},\F)$ relative to $\nu_{2^n}$. Put
${\cal N}_1:= {\cal N}\cap W^{\mathrm{vs}}_n$ and ${\cal N}_2:={\langle n_{l_X}\rangle}_{X\in\Delta_{n-1}}$, where $X$ ranges in the set of $(n-1)$-elements of $\Delta$, $l_X$ is the line of $\Delta_n$ corresponding
to $X$ (see (\ref{lines}) of Subsection \ref{1.2}) and $n_{l_X}$ is the nucleus of the conic $\varepsilon_n^{\mathrm{vs}}(l_X)$. Clearly, $\mathcal{N}_1\supseteq \mathcal{N}_2$ and both these subspaces are stabilized by the group $G = \mathrm{Spin}(2n+1,\F)$. We can also define two mappings $\varepsilon_n^{\mathrm{vs}}/{\cal N}_1$ and $\varepsilon_n^{\mathrm{vs}}/{\cal N}_2$ from $\Delta_n$ to the set of $1$-dimensional linear subspaces of $W_n^{\mathrm{vs}}/{\cal N}_1$ and $W_n^{\mathrm{vs}}/{\cal N}_2$ respectively and a mapping $\iota_{n-1}:\Delta_{n-1}\rightarrow{\cal N}_2$, as follows: $(\varepsilon_n^{\mathrm{vs}}/{\cal N}_1)(X) := \langle \varepsilon_n^{\mathrm{vs}}(X), {\cal N}_1\rangle/{\cal N}_1$ and $(\varepsilon_n^{\mathrm{vs}}/{\cal N}_2)(X) := \langle \varepsilon_n^{\mathrm{vs}}(X), {\cal N}_2\rangle/{\cal N}_2$ for every point $X$ of $\Delta_n$ and
$\iota_{n-1}(X) := n_{l_X}$ for every $(n-1)$-element $X$ of $\Delta$.

The next lemma, to be proved in Section \ref{sec3}, allows us to define one more $G$-invariant subspace of $W^{\mathrm{vs}}_n$.

\begin{Lemma}\label{main lemma n = 2}
Let $n = 2$. Then $\iota_1\cong\tilde{\varepsilon}_1$ $(\cong \varepsilon_1^{\mathrm{gr}})$.
\end{Lemma}
In other words, if $n = 2$ then $\iota_1(\Delta_1)$ is a copy of the quadric $\Delta_1\cong Q(4,\F)$ in $\PG({\cal N}_2)\cong \PG(4,\F)$. One of the points of $\PG({\cal N}_2)$ is the nucleus of the quadric $\iota_1(\Delta_1)$.

Let now $n > 2$. Given a $k$-element $X$ of $\Delta$ with $k \leq n-2$, let $\mathrm{Res}^+(X)$ be the {\em upper residue} of $X$ in $\Delta$, formed by the elements of $\Delta$ of type $i > k$ incident with $X$. Then $\mathrm{Res}^+(X)$ is a building of type $B_{n-k}$ with $\{k+1,..., n\}$ as the set of types. We can define the $n$-{\em grassmannian} $\mathrm{Res}^+_n(X)$ of $\mathrm{Res}^+(X)$ by taking the $n$-elements of $\mathrm{Res}^+(X)$ as points of $\mathrm{Res}^+_n(X)$ and the lines of $\Delta_n$ contained in $\mathrm{Res}^+(X)$ as lines of $\mathrm{Res}^+_n(X)$. When $k = n-2$ we denote $\mathrm{Res}_n^+(X)$ by the symbol $Q_X$ and we call it a {\em quad} of $\Delta_n$. We have $\mathrm{dim}(\langle \varepsilon^{\mathrm{spin}}_n(Q_X)\rangle) = 4$, namely $\varepsilon^{\mathrm{spin}}_n$ embeds $Q_X$ as a copy of the symplectic generalized quadrangle $W(3,\F)$ in the 4-space $\langle \varepsilon^{\mathrm{spin}}_n(Q_X)\rangle$. By Lemma \ref{main lemma n = 2}, $\iota_{n-1}(Q_X)$ is a copy
of $Q(4,\F)$ in the 5-dimensional subspace $\langle\iota_{n-1}(Q_X)\rangle$ of $\PG(W^{\mathrm{vs}}_n)$. We denote by $n_{Q_X}$ the nucleus of the quadric $\iota_{n-1}(Q_X)$ in $\langle\iota_{n-1}(Q_X)\rangle$ and we put ${\cal N}_3 = \langle n_{Q_X}\rangle_{X\in\Delta_{n-2}}$. Clearly, ${\cal N}_3$ is a subspace of ${\cal N}_2$ and it is stabilized by $G$.

We denote by $\iota_{n-2}:\Delta_{n-2}\rightarrow{\cal N}_3$ the mapping which maps every $X\in \Delta_{n-2}$ onto $n_{Q_X}$. We also denote by  $\varepsilon_n^{\mathrm{vs}}/{\cal N}_3$ the function which maps $X$ of $\Delta_n$ onto $\langle \varepsilon_n^{\mathrm{vs}}(X), {\cal N}_3\rangle/{\cal N}_3$ for every point $X$ of $\Delta_n$. We extend this latter notation to the case of $n = 2$ by stating that in that case ${\cal N}_3$ is the nucleus of $\iota_1(\Delta_1)$.

\begin{theorem}\label{main theorem 2}
Let $\F$ be a perfect field of characteristic $2$. Then:\\
(1) The mappings $\varepsilon_n^{\mathrm{vs}}/{\cal N}_1$ and $\varepsilon_n^{\mathrm{vs}}/{\cal N}_2$ are projective embeddings of $\Delta_n$ while $\varepsilon_n^{\mathrm{vs}}/{\cal N}_3$ is a veronesean embedding. \\
(2) $\varepsilon_n^{\mathrm{vs}}/{\cal N}_1\cong \varepsilon^{\mathrm{spin}}_n$. \\
(3) The embeddings $\varepsilon_n^{\mathrm{sp}}$ and $\varepsilon_n^{\mathrm{gr}}$ are (possibly improper) quotients of $\varepsilon_n^{\mathrm{vs}}/{\cal N}_2$ and $\varepsilon_n^{\mathrm{vs}}/{\cal N}_3$ respectively.\\
(4) The mapping $\iota_{n-1}$ is a projective embedding of the $(n-1)$-grassmannian $\Delta_{n-1}$ of $\Delta$ in ${\cal N}_2$. \\
(5) Let $n > 2$. Then $\iota_{n-2}$ is a projective embedding of the $(n-2)$-grassmannian $\Delta_{n-2}$ of $\Delta$ in ${\cal N}_3$.
\end{theorem}
We shall prove this theorem in Section \ref{sec3}. Perhaps some of the claims gathered in Theorem \ref{main theorem 2} also hold when $\F$ is non-perfect. In particular, we make no use of the hypothesis that $\F$ is perfect in the proof of (4). On the other hand, if $\F$ is non-perfect then $\varepsilon_n^{\mathrm{vs}}/{\cal N}_1$ and $\varepsilon_n^{\mathrm{vs}}/{\cal N}_2$ are lax embeddings. Thus, (2) and that part of (3) that deals with $\varepsilon_n^{\mathrm{vs}}/{\cal N}_2$ are false when $\F$ is non-perfect.

The next theorem, which also will be proved in Section \ref{sec3}, contains a complete description of the module $W^{\mathrm{vs}}_n \cong V(2\lambda_n)$ when $n \leq 3$. We will exploit it to prove some of the claims of Theorem \ref{main theorem 2}.

\begin{theorem}\label{main theorem 4new}
Let $\F$ be a perfect field of characteristic $2$. \\
(1) Let $n = 2$. Then $\varepsilon_2^{\mathrm{vs}}/{\cal N}_2\cong\varepsilon_2^{\mathrm{sp}}$ and $\varepsilon_2^{\mathrm{vs}}/{\cal N}_3\cong\varepsilon_2^{\mathrm{gr}}$. The $G$-module
$W^{\mathrm{vs}}_2$ admits the following composition series: $W^{\mathrm{vs}}_2\supset{\cal N}_1\supset{\cal N}_2\supset{\cal N}_3\supset 0$. We have $\mathrm{dim}(W^{\mathrm{vs}}_2/{\cal N}_1) = 4$, $\mathrm{dim}({\cal N}_1/{\cal N}_2) = 1$, $\mathrm{dim}({\cal N}_2/{\cal N}_3) = 4$ and $\mathrm{dim}({\cal N}_3) = 1$. \\
(2) Let $n = 3$. Then $\varepsilon_3^{\mathrm{vs}}/{\cal N}_2\cong\varepsilon_3^{\mathrm{sp}}$, $\varepsilon_3^{\mathrm{vs}}/{\cal N}_3\cong \varepsilon_3^{\mathrm{gr}}$, $\iota_2\cong\tilde{\varepsilon}_2$ and $\iota_1\cong\tilde{\varepsilon}_1$ $(\cong \varepsilon^{\mathrm{gr}}_1)$.
The $G$-module $W^{\mathrm{vs}}_3$ admits a composition series as follows: $W^{\mathrm{vs}}_3\supset{\cal N}_1\supset{\cal N}_2\supset{\cal N}_3\supset K_0 \supset 0$ where $\mathrm{dim}(K_0) = 1$. Moreover, $\mathrm{dim}(W^{\mathrm{vs}}_3/{\cal N}_1) = 8$, $\mathrm{dim}({\cal N}_1/{\cal N}_2) = 6$, $\mathrm{dim}({\cal N}_2/{\cal N}_3) = 14$ and $\mathrm{dim}({\cal N}_3/K_0) = 6$.
\end{theorem}
In Section 4 we address the problem of the universality of $\tilde{\varepsilon}_n$. We focus on the case of $n = 2$, proving the following two theorems.

\begin{theorem}\label{main theorem univ1}
Let $n = 2$. Let $\F$ be a finite field of odd order $|\F| > 3$. Then $\tilde{\varepsilon}_2$ is relatively universal.
\end{theorem}

\begin{theorem}\label{main theorem univ2}
Let $n = 2$ and $\mathrm{char}(\F) = 2$. Then $\tilde{\varepsilon}_2$ is not relatively universal.
\end{theorem}
We are presently unable to say so much on the case of $n > 2$. The following is the only result that we can offer for that case.

\begin{theorem}\label{main theorem univ3}
Let $\F$ be a perfect field of characteristic $2$. Then $\tilde{\varepsilon}_n$ is not relatively universal, for any $n \geq 2$.
\end{theorem}

\subsection{Conjectures and problems}
 \noindent
{\bf Conjectures.} (1) Let $\F$ be a perfect field of characteristic 2 and let $n > 3$. Then $\tilde{\varepsilon}_n/{\cal N}_2\cong \varepsilon_n^{\mathrm{sp}}$, $\tilde{\varepsilon}_n/{\cal N}_3\cong \varepsilon_n^{\mathrm{gr}}$, $\iota_{n-1}\cong\tilde{\varepsilon}_{n-1}$ and $\iota_{n-2}\cong\tilde{\varepsilon}_{n-2}$. \\
(2) We have $\tilde{\varepsilon}_n/{\cal N}_3\cong \varepsilon_n^{\mathrm{gr}}$, $\iota_{n-1}\cong\tilde{\varepsilon}_{n-1}$ and $\iota_{n-2}\cong\tilde{\varepsilon}_{n-2}$ even if $\F$ is non-perfect.\\
(3) The hypothesis that $\F$ is perfect is superfluous in Theorem \ref{main theorem univ3}.\\
(4) If $\mathrm{char}(\F)\neq 2$ and $\F \neq \F_3$ then $\tilde{\varepsilon}_n$ is relatively universal, for any $n\geq 2$.\\

\noindent
{\bf Problems.} (1) Is the assumption $|\F| > 3$ really necessary in Theorem \ref{main theorem univ1}?\\
(2) Does $\Delta_n$ admit the absolutely universal veronesean embedding?

\section{Proof of Theorem \ref{main theorem 1}}\label{sec2}

In this section we will freely use various notions and results from the theory of Chevalley groups and their Lie algebras. We mainly rely on Steinberg \cite{Stein}, Carter \cite{Carter} and Humphreys \cite{Humphreys}, \cite{Humphreys2} for this matter.

We recall that Chevalley groups are firstly defined over the complex field $\mathbb{C}.$ Given a (simple) Lie algebra $\mathfrak{L}_{\mathbb{C}}$ over $\mathbb{C}$ and an irreducible $\mathfrak{L}_{{\mathbb{C}}}$-module $V_{\mathbb{C}}(\lambda)$ with highest weight $\lambda$, the subgroup $G_{\mathbb{C}}$ of $\mathrm{SL}(V_{\mathbb{C}}(\lambda))$ generated by the exponential maps $e^{X_{\alpha}t}$ (for $\alpha$ a root) is the complex Chevalley group  associated to $\mathfrak{L}_{\mathbb{C}}$ (and $\lambda$). In order to replace $\mathbb{C}$ with an arbitrary field $\F$, a suitable lattice $L$ is chosen in $V_{\mathbb{C}}(\lambda).$
The analogue $G_{\F}$ of $G_{\mathbb{C}}$ is defined in $V_{\F}(\lambda):=\F\otimes_{\mathbb{Z}}L$ in the same way as $G_{\mathbb{C}}$ in $V_{\mathbb{C}}(\lambda).$ (See Steinberg~\cite{Stein} or Carter~\cite{Carter} for more details). Chevalley groups can also be obtained in a different way, as $\F$-rational subgroups of (almost simple) algebraic groups defined over the algebraic closure $\overline{\F}$ of $\F$ (see e.g. Humphreys~\cite{Humphreys2}.) We will not make use of this latter perspective in this paper except in one occasion in Subsection~\ref{subsection} (and with $\overline{\F}=\F=\mathbb{C}$).

\subsection{Preliminaries}

\subsubsection{A result on quadric veronesean maps}

We will firstly prove an analog of Theorem \ref{main theorem 1} for the veronesean embedding $\nu_m:\PG(m-1,\F)\rightarrow V({{m+1}\choose 2},\F)$. The embedding $\nu_m$ is $\mathrm{SL}(m,\F)$-homogeneous. Thus $V({{m+1}\choose 2},\F)$ acquires the structure of an $\mathrm{SL}(m,\F)$-module. Let ${\cal W}_\F$ denote this $\mathrm{SL}(m,\F)$-module.

Let $\omega_1, \omega_2,..., \omega_{m-1}$ be the fundamental dominant weights for the root system of type $A_{m-1}$, $m\geq 3$. The nodes of the diagram $A_{m-1}$ are numbered in the usual way:

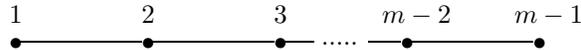
\begin{picture}(310,36)(0,0)
\put(20,8){$\bullet$}
\put(23,11){\line(1,0){47}}
\put(70,8){$\bullet$}
\put(73,11){\line(1,0){47}}
\put(120,8){$\bullet$}
\put(123,11){\line(1,0){12}}
\put(138,10){$.....$}
\put(156,11){\line(1,0){12}}
\put(168,8){$\bullet$}
\put(171,11){\line(1,0){47}}
\put(218,8){$\bullet$}
\put(20,18){1}
\put(70,18){2}
\put(120,18){3}
\put(161,18){$m-2$}
\put(211,18){$m-1$}
\end{picture}

\noindent
Put $\omega = 2\cdot \omega_1$ and let $V_\F(\omega)$ be the Weyl module for $\mathrm{SL}(m,\F)$ with $\omega$ as the highest weight. We put $\F$ as a subscript in the symbol $V_\F(\omega)$ in order to keep a record of the field $\F$ in our notation. We shall adopt this expedient for nearly all symbols in this section.

\begin{theo}\label{projective veronesean}
The $\mathrm{SL}(m,\F)$-modules ${\cal W}_\F$ and $V_\F(\omega)$ are isomorphic for all $m \geq 3$ and all fields $\F$.
\end{theo}
\pr We firstly assume that $\F = \mathbb{C}$. Let $\mathfrak{L}_{\mathbb{C}} = \mathfrak{sl}(m,\mathbb{C}) = \mathfrak{H}\oplus(\oplus_{\alpha\in \Phi}X_{\alpha})$ be the Lie algebra of $\mathrm{SL}(m,\mathbb{C})$, where $\mathfrak{H}$ is the Cartan subalgebra, $\Phi$ is the set of all roots and $X_{\alpha}$ is the $1$-dimensional subalgebra of $\mathfrak{L}_{\mathbb{C}}$ corresponding to the root $\alpha$. Let $Z_{\mathbb{C}}(\omega)$ be the cyclic $\mathfrak{L}_{\mathbb{C}}$-module associated to $\omega$ and $v_0$ a highest weight vector of $Z_{\C}(\omega)$. We recall that $Z_{\C}(\omega)$ admits a unique maximal proper submodule $J_{\C}(\omega)$ and $V_{\C}(\omega) = Z_{\C}(\omega)/J_{\C}(\omega)$.

As ${\cal W}_{\C}$ is an $\mathrm{SL}(m,\C)$-module the algebra $\mathfrak{L}_{\C}$ also acts on ${\cal W}_{\C}$. Explicitly, given a basis $\{e_1,..., e_m\}$ of $V(m,\C)$ the vectors of ${\cal W}_{\C}$ can be regarded as linear combinations of the vectors $e_i\otimes e_i$ for $i = 1,...,m$ and $e_i\otimes e_j+e_j\otimes e_i$ for $1\leq i < j \leq m$. Thus, the set
\[\{e_i\otimes e_i\}_{i=1}^m\cup\{e_i\otimes e_j+e_j\otimes e_i\}_{1\leq i < j \leq m}\]
can be taken as a basis of ${\cal W}_\C$. Identifying $\mathfrak{L}_{\mathbb{C}} = \mathfrak{sl}(m,\mathbb{C})$ with the Lie algebra of endomorphisms of $V(m,\C)$ with null trace, if $a\in\mathfrak{L}_{\C}$ then
\[\begin{array}{l}
a(e_i\otimes e_i) = a(e_i)\otimes e_i + e_i\otimes a(e_i),\\
a(e_i\otimes e_j +  e_j\otimes e_i) = a(e_i)\otimes e_j + e_i\otimes a(e_j) + a(e_j)\otimes e_i + e_j\otimes a(e_i).
\end{array}\]
We can assume to have chosen the decomposition $\mathfrak{L}_{\mathbb{C}} = \mathfrak{H}\oplus(\oplus_{\alpha\in \Phi}X_{\alpha})$ so that for $i = 1, 2,..., m-1$ the $1$-dimensional summand $X_{\alpha_i}$ corresponding to the simple root $\alpha_i$ is the set of endomorphisms $a$ of $V(m,\C)$ such that $a(e_{i+1}) \in \langle e_i\rangle$ and $a(e_j) = 0$ for $j \neq i+1$.

Let $\bar{v}_0 := \nu_m(e_1)$. It is straightforward to check that $h(\bar{v}_0) = \omega(h)\cdot\bar{v}_0$ for every $h\in\mathfrak{H}$ and $X_\alpha(\bar{v}_0) = 0$ for every positive root $\alpha$. Moreover, the $\mathfrak{L}_{\C}$-orbit of $\bar{v}_0$ spans ${\cal W}_{\C}$. Hence there exists a surjective homomorphism
of $\mathfrak{L}_{\C}$-modules $f_{\C}$ from $Z_{\C}(\omega)$ to ${\cal W}_{\C}$ mapping $v_0$ onto $\bar{v}_0$. The kernel $\mathrm{ker}(f_{\C})$ of $f_{\C}$ is a proper submodule of $Z_{\C}(\omega)$. Hence $\mathrm{ker}(f_{\C})\subseteq J_{\C}(\omega)$. It follows that $\mathrm{dim}({\cal W}_{\C}) \geq \mathrm{dim}(Z_{\C}(\omega)/J_{\C}(\omega))$, namely $\mathrm{dim}({\cal W}_{\C}) \geq \mathrm{dim}(V_{\C}(\omega))$. Clearly, ${\cal W}_{\C} \cong V_{\C}(\omega)$ if and only if $\mathrm{dim}({\cal W}_{\C}) = \mathrm{dim}(V_{\C}(\omega))$. However $\mathrm{dim}(V_{\C}(\omega)) = {{m+1}\choose 2}$, as one can see by using the Weyl dimension formula. Hence $\mathrm{dim}(V_{\C}(\omega)) = \mathrm{dim}({\cal W}_{\C})$, namely ${\cal W}_{\C} \cong V_{\C}(\omega)$. Accordingly, $f_{\C}$ splits as $f_{\C} = \varphi_{\C}\circ\pi_{\C}$ where $\pi_{\C}$ is the canonical projection of $Z_{\C}(\omega)$ onto $V_{\C}(\omega)$ and $\varphi_{\C}$ is the isomorphism from $V_{\C}(\omega)$ to ${\cal W}_{\C}$ mapping $v_0 + J_{\C}(\omega)$ onto $\bar{v}_0$.

Before replacing $\C$ with an arbitrary field $\F$ we must turn from $\C$ to the ring of integers $\mathbb{Z}.$ To this aim, let $\mathfrak{U}_{\C}$ be the enveloping (associative) algebra of $\mathfrak{L}_{\C}$, let $\{\alpha_1,..., \alpha_{m-1}\}$ be the set of simple roots and choose a Chevalley basis $\{h_1,..., h_{m-1}, x_1,...,x_{m-1}, y_1,..., y_{m-1}\}$ of $\mathfrak{L}_{\C}$, where $h_i\in \mathfrak{H}$, $x_i\in X_{\alpha_i}$ and $y_i\in X_{-\alpha_i}$ for $i = 1,...,m-1$. Regarding $\mathfrak{U}_{\C}$ as a ring, let $\mathfrak{U}_{\mathbb{Z}}$ be the subring of $\mathfrak{U}_{\C}$ generated by the elements of the following form
\[\frac{x_1^{r_1}}{r_1!}\cdot...\cdot\frac{x_{m-1}^{r_{m-1}}}{r_{m-1}!}\cdot{{h_1}\choose{s_1}}\cdot...\cdot{{h_{m-1}}\choose{s_{m-1}}}\cdot...\cdot\frac{y_1^{t_1}}{t_1!}\cdot...\cdot\frac{y_{m-1}^{t_{m-1}}}{t_{m-1}!}\]
for nonnegative integers $r_1,..., r_{m-1}, s_1,..., s_{m-1}, t_1,..., t_{m-1}$. We recall that
\[{{h}\choose{s}} := \frac{h(h-\iota)...(h-(s-1)\iota)}{s!}\]
where $\iota$ is the identity element of $\mathfrak{U}_{\C}$. Let $Z_{\mathbb{Z}}(\omega)$, $V_{\mathbb{Z}}(\omega)$ and ${\cal W}_{\mathbb{Z}}$ be the $\mathfrak{U}_{\mathbb{Z}}$-orbits of $v_0$, $v_0+J_{\C}(\omega)$ and $\bar{v}_0$ respectively. Then $Z_{\mathbb{Z}}(\omega)$, $V_{\mathbb{Z}}(\omega)$ and ${\cal W}_{\mathbb{Z}}$ are left $\mathfrak{U}_{\mathbb{Z}}$-modules, $f_{\C}$ induces a homomorphism $f_{\mathbb{Z}}$ from $Z_{\mathbb{Z}}(\omega)$ to ${\cal W}_{\mathbb{Z}}$ and $\varphi_{\C}$ induces an isomorphism from $V_{\mathbb{Z}}(\omega)$ to ${\cal W}_{\mathbb{Z}}$. Also, $\mathrm{ker}(f_{\mathbb{Z}}) = J_{\mathbb{Z}}(\omega) := Z_{\mathbb{Z}}(\omega)\cap J_{\C}(\omega)$, whence $V_{\mathbb{Z}}(\omega) = Z_{\mathbb{Z}}(\omega)/J_{\mathbb{Z}}(\omega)$. Moreover $V_{\mathbb{Z}}(\omega)$ is a lattice in $V_{\C}(\omega)$. Let $B$ be a basis of the lattice $V_{\mathbb{Z}}(\omega)$. Thus $B$ is a basis of the $\mathbb{Z}$-module $V_{\mathbb{Z}}(\omega)$ as well as a basis of the $\C$-vector space $V_{\C}(\omega)$. Since $\varphi_{\C}$ is an isomorphism, the set $\overline{B} := \varphi_{\mathbb{Z}}(B)$ is a basis of the $\mathbb{Z}$-module ${\cal W}_{\mathbb{Z}}$ as well as of the $\C$-vector space ${\cal W}_{\C}$.

Let now $\F$ be an arbitrary field. Then $V_\F(\omega) = \F\otimes_{\mathbb{Z}}V_{\mathbb{Z}}(\omega)$, where integers are taken modulo $p = \mathrm{char}(\F)$ if $\mathrm{char}(\F) > 0$. Note that $B$ is a basis of the $\F$-vector space $\F\otimes_{\mathbb{Z}}V_{\mathbb{Z}}(\omega)$. Keeping the symbol $\bar{v}_0$ to denote the vector of ${\cal W}_\F$ defined in the same way as $\bar{v}_0$ in ${\cal W}_\C$, the $\mathfrak{U}_{\mathbb{Z}}$-orbit of $\bar{v}_0$ spans ${\cal W}_\F$. Hence ${\cal W}_{\F} = \F\otimes_{\mathbb{Z}}{\cal W}_{\mathbb{Z}}$, the set $\overline{B}$ is a basis of ${\cal W}_\F$ and $\varphi_{\mathbb{Z}}$ uniquely determines an isomorphism $\varphi_\F:V_\F(\omega)\rightarrow{\cal W}_{\F}$. \eop

\bigskip

\noindent
{\bf Remark.} When $\mathrm{char}(\F) = 2$ the module ${\cal W}_{\F}$ is reducible. Indeed it admits an ${{m}\choose 2}$-dimensional submodule, spanned by the vectors $e_i\otimes e_j+e_j\otimes e_i$ for $i < j$. This submodule corresponds to the submodule of $V_\F(\omega)$ spanned by the $SL(m,\F)$-orbit of $X_{-\alpha_1}(\bar{v}_0) = X_{-\alpha_1}(\nu_m(e_1)) = X_{-\alpha_1}(e_1\otimes e_1) = \{t(e_2\otimes e_1 + e_1\otimes e_2)\}_{t\in \F}$.

\subsubsection{Notation}

Let $m = 2^n$ and let $\omega_1,..., \omega_{2^n-1}$ be the fundamental dominant weights of the root system of type $A_{2^n-1}$. Let $\lambda_1,..., \lambda_n$ be the fundamental dominant weights for the root system of type $B_n$, as in the introduction of this paper. The vector space $V(2^n,\F)$ supports both the Weyl module $V_\F(\omega_1)$ for $\mathrm{SL}(2^n,\F)$ and the spin module $V_\F(\lambda_n)$. Clearly, $V_\F(2\lambda_n)$ is the module denoted by $V(2\lambda_n)$ in the introduction of this paper.

We denote the group $\mathrm{Spin}(2n+1,\F)$ by $G_\F$ and we put $L_\F := \mathrm{SL}(2^n,\F)$, for short. We will use the symbol $W^{\mathrm{vs}}_\F$ to denote the $G_\F$-module $\langle \varepsilon^{\mathrm{vs}}_n(V_\F(\lambda_n))\rangle$ and we denote the $L_\F$-module $\langle \nu_{2^n}(V_\F(\omega_1)\rangle$ by ${\cal W}_\F$, consistently with the notation of the previous subsection. Clearly, ${\cal W}_\F$ is also a (reducible) $G_\F$-module. We recall that ${\cal W}_\F \cong V_\F(2\omega_1)$ by Theorem \ref{projective veronesean}.

\subsubsection{A few results on $\mathfrak{o}(2n+1,\C)$ and $\mathfrak{sl}(2^n,\C)$}\label{subsection}

The group $G_\F$ acts faithfully on $V(2^n,\F)$ as a subgroup of $L_\F$. So, for the rest of this section we regard $G_\F$ as a subgroup of $L_\F$.

\begin{lemma}\label{tori}
For any field $\F$, every maximal split torus of $G_\F$ is contained in a maximal split torus of $L_\F$.
\end{lemma}
\pr It suffices to prove that every maximal split torus of $G_\F$ stabilizes all $1$-dimensional subspaces $\langle e_1\rangle,..., \langle e_{2^n}\rangle$ for a suitable basis $\{e_i\}_{i=1}^{2^n}$ of $V(2^n,\F)$, where $\langle e_1\rangle,..., \langle e_{2^n}\rangle$ are the images under $\varepsilon^{\mathrm{spin}}$ of the $n$-elements of an apartment of the $B_n$-building $\Delta$.

If $n = 2$ then $G_\F$ is the symplectic group $\mathrm{Sp}(4,\F)$ in its natural action on $V(4,\F)$. In this case the claim is obvious. Let $n > 2$ and let $T$ be a maximal split torus of $G_\F$. Then $T$ is the stabilizer of two opposite chambers $C_1$ and $C_2$ of $\Delta$. Equivalently, $T$ stabilizes all chambers of the unique apartment of $\Delta$ contaning $C_1$ and $C_2$. Let $p_1$ and $p_2$ be the elements of type $1$ of $C_1$ and $C_2$ respectively. Let $C'_1$ be the projection of $C_2$ into the residue $\Delta' := \mathrm{Res}_\Delta(p_1)$ of $p_1$ (Tits \cite[Chapter 3]{Tits}). We recall that $\Delta'$ is a building of type $B_{n-1}$. The flags $C_1\setminus\{p_1\}$ and $C'_1\setminus\{p_1\}$ are opposite chambers of $\Delta'$ and $T$ stabilizes both $C_1$ and $C'_1$. On the other hand, the spin embedding $\varepsilon_n^{\mathrm{spin}}$ of the dual polar space $\Delta_n$ in $V(2^n,\F)$ induces the spin embedding $\varepsilon_{n-1}^{\mathrm{spin}}$ on the dual polar space associated to $\Delta'$. By the inductive hypothesis, $T$ fixes all $1$-dimensional subspaces $\langle e_{1,i}\rangle$ for a basis $\{e_{1,i}\}_{i=1}^{2^{n-1}}$ of $\langle\varepsilon_n^{\mathrm{spin}}(\Delta_n(p_1))\rangle$, where $\Delta_n(p_1)$ stands for the set of points of $\Delta'$, namely the set of $n$-elements of $\Delta$ incident to $p_1$. Moreover, taking $\{2,3,..., n\}$ as the set of types of $\Delta'$, the subspaces $\langle e_{1,1}\rangle,..., \langle e_{1,2^{n-1}}\rangle$ correspond to the $n$-elements of the apartment of $\Delta'$ containing $C_1\setminus\{p_1\}$ and $C'_1\setminus \{p_1\}$. Similarly, $T$ fixes all $1$-dimensional subspaces $\langle e_{2,i}\rangle$ for a basis $\{e_{2,i}\}_{i=1}^{2^{n-1}}$ of $\langle\varepsilon_n^{\mathrm{spin}}(\Delta_n(p_2))\rangle$ and the subspaces $\langle e_{2,1}\rangle,..., \langle e_{2,2^{n-1}}\rangle$ correspond to the $n$-elements of the apartment of $\mathrm{Res}_\Delta(p_2)$ contaning $C_2\setminus \{p_2\}$ and $C'_2\setminus\{p_2\}$, where $C'_2$ is the projection of $C_1$ into $\mathrm{Res}_\Delta(p_2)$.

The $1$-subspaces $\langle e_{1,1}\rangle,..., \langle e_{1,2^{n-1}}\rangle, \langle e_{2,1}\rangle,..., \langle e_{2,2^{n-1}}\rangle$ correspond to the $n$-elements of the apartment of $\Delta$ containing $C_1$ and $C_2$. Moreover the dual polar space $\Delta_n$ is generated by $\Delta_n(p_1)\cup\Delta_n(p_2)$ (see Blok and Brouwer \cite{BB}, also Cooperstein and Shult \cite{CS}). On the other hand, $V(2^n,\F)$ is spanned by $\langle\varepsilon_n^{\mathrm{spin}}(\Delta_n)\rangle$. Therefore $\{e_{1,i}\}_{i=1}^{2^{n-1}}\cup\{e_{2,i}\}_{i=1}^{2^{n-1}}$ is a basis of $V(2^n,\F)$.  \eop

\bigskip

Let now $\F = \C$. We denote the Lie algebras of $L_\C$ and $G_\C$ by the symbols $\mathfrak{L}_\C$ and $\mathfrak{L}^\circ_\C$ respectively. Explicitly, $\mathfrak{L}_\C = \mathfrak{sl}(2^n,\C)$ and $\mathfrak{L}^\circ = \mathfrak{o}(2n+1,\C)$. We use the symbol $\mathfrak{H}_\C$ to denote a Cartan subalgebra of $\mathfrak{L}_\C$ and $\mathfrak{H}^\circ_\C$ for a Cartan subalgebra of $\mathfrak{L}^\circ_\C$.

We have $G_\C \subset L_\C$. It is well known that $G_\C$ is a closed subgroup of $L_\C$ in the Zariski topology (see e.g. Cassels \cite[Chapter 10]{Cass}). The Lie algebra $\mathfrak{L}^\circ_\C$ is supported by the tangent space of $G_\C$ at the identity element, which is a subspace of the tangent space of $L_\C$ at the identity element. Thus the inclusion $G_\C\subset L_\C$ naturally induces an embedding $\iota:\mathfrak{L}^\circ_\C\rightarrow\mathfrak{L}_\C$. The Cartan subalgebras of $\mathfrak{L}^\circ_\C$ bijectively correspond to the maximal tori of $G_\C$. By Lemma \ref{tori}, every Cartan subalgebra $\mathfrak{H}^\circ_\C$ of $\mathfrak{L}^\circ_\C$ is mapped under $\iota$ into a Cartan subalgebra $\mathfrak{H}_\C$ of $\mathfrak{L}_\C$ such that for a suitable basis $B$ of $V(2^n,\C)$ the algebras $\mathfrak{H}^\circ_\C$ and $\mathfrak{H}_\C$ stabilize the $1$-linear subspaces spanned by the vectors of $B$. With this choice of $\mathfrak{H}^\circ_\C$ and $\mathfrak{H}_\C$ we have the following:

\begin{lemma}\label{pesi}
The weight $\lambda_n$ is the restriction of $\omega_1$ to $\mathfrak{H}_{\C}^\circ$. Similarly, $2\cdot\lambda_n$ is the restriction of $2\cdot\omega_1$ to $\mathfrak{H}_{\C}^\circ$.
\end{lemma}
\pr Let $u_0$ be a highest weight vector of $V_\C(\lambda_n)$. We may assume that, given two opposite chambers $C_1$ and $C_2$ of $\Delta$ as in the proof of Lemma \ref{tori}, the vector $u_0$ spans the image under $\varepsilon_n^{\mathrm{spin}}$ of the element of $C_1$ of type $n$ and that $\mathfrak{H}_{\C}^\circ$ is associated to the maximal split torus of $G_\C$ stabilizing all chambers of the apartment $\cal A$ of $\Delta$ containing $C_1$ and $C_2$. We can also assume that $\mathfrak{H}_\C$ is associated to the maximal split torus of $L_\C$ stabilizing the images under $\varepsilon^{\mathrm{spin}}_n$ of the $n$-elements of $\cal A$ (see the proof of Lemma \ref{tori}). Hence $\mathfrak{H}_{\C}^\circ\subseteq \mathfrak{H}_{\C}$ and $u_0$ can be taken as a highest weight vector of $V_\C(\omega_1)$ too. Since $\mathfrak{H}_{\C}^\circ\subseteq \mathfrak{H}_{\C}$, we have $h(u_0)=\lambda_n(h)\cdot u_0 = \omega_1(h)\cdot u_0$ for any $h\in \mathfrak{H}^\circ_\C$. This proves the first claim of the Lemma. The second claim immediately follows from the first one. \eop

\subsection{Proof of Theorem \ref{main theorem 1}}

\subsubsection{The case $\F = \C$}

Modulo a suitable choice of a $BN$-pair in $G_\C$ and $L_\C$, any of the vectors representing a point of $\Delta_n$ can be chosen as a highest weight vector of $V_\C(\lambda_n)$ and any non-zero vector of $V(2^n,\C)$ can be chosen as a highest weight vector for $V_\C(\omega_1)$. Thus, we may assume to have chosen the same vector $u_0$ of $V(2^n,\C)$ as a highest weight vector of either $V_\C(\lambda_n)$ and $V_\C(\omega_1)$ as in the proof of Lemma \ref{pesi}. The quadric veronesean map $\nu_{2^n}$ is defined modulo the choice of a basis $B$ of $V(2^n,\C)$. We may assume to have chosen $B$ in such a way that $u_0$ is the first element of $B$. So, $\bar{v}_0 := \nu_{2^n}(u_0) = u_0\otimes u_0$.

Let $v_0$ be a highest weight vector of $V_\C(2\omega_1)$ and $v^\circ_0$ a highest weight vector of $V_\C(2\lambda_n)$. As in the proof of Theorem \ref{projective veronesean}, let $f_\C$ be the homomorphism from $Z_\C(2\omega_1)$ to ${\cal W}_\C$ mapping $v_0$ onto $\bar{v}_0$ and let $\varphi_\C$ be the isomorphism from $V_\C(2\omega_1)$ to ${\cal W}_\C$ induced by $f_\C$.

Denote by $\Phi$ and $\Phi^\circ$ the root systems of type $A_{2^n-1}$ and $B_n$ respectively, let $\Phi_+$ and $\Phi^\circ_+$ be the set of positive roots in $\Phi$ and $\Phi^\circ$. Then
\[\mathfrak{L}_\C  =  \mathfrak{H}_\C \oplus (\oplus_{\alpha\in \Phi}X_\alpha), \hspace{8 mm} \mathfrak{L}^\circ_\C  =  \mathfrak{H}^\circ_\C \oplus (\oplus_{\alpha\in \Phi^\circ}X^\circ_\alpha),\]
where $X_\alpha$ and $X^\circ_\alpha$ are the 1-dimensional subalgebras of $\mathfrak{L}_\C$ and $\mathfrak{L}^\circ_\C$ associated to the roots.

As $u_0$ is a highest weight vector for both $V_\C(\omega_1)$ and $V_\C(\lambda_n)$, we have $X_{\alpha}(u_0) = X^\circ_\beta(u_0) = 0$ for any $\alpha\in \Phi_+$ and any $\beta\in\Phi^\circ_+$. As $\bar{v}_0 = u_0\otimes u_0$, the above implies that $X_\alpha(\bar{v}_0) = X^\circ_\beta(\bar{v}_0) = 0$ for any $\alpha\in \Phi_+$ and any $\beta\in\Phi^\circ_{+}$. Moreover, $h(\bar{v}_0) = 2\omega_1(h)\cdot \bar{v}_0 = 2\lambda_n(h)\cdot \bar{v}_0$ for every $h\in\mathfrak{H}_\C^\circ$ by Corollary \ref{pesi}. Therefore the submodule of ${\cal W}_\C$ spanned by the orbit of $\bar{v}_0$ under the action of $\mathfrak{L}^\circ_\C$ is a homomorphic image of the cyclic $\mathfrak{L}^\circ_\C$-module $Z_\C(2\lambda_n)$. However the submodule of ${\cal W}_\C$ spanned by the $\mathfrak{L}^\circ_\C$-orbit of $\bar{v}_0$ is just $W^{\mathrm{vs}}_\C$. Thus, there exists a surjective homomorphism $f^\circ_\C$ from $Z_\C(2\lambda_n)$ to $W^{\mathrm{vs}}_\C$. This homomorphisms splits as $f^\circ_\C = {\varphi^\circ_\C}\circ g^\circ_\C$, where $\varphi^\circ_\C$ is the isomorphism from $\varphi_\C^{-1}(W^{\mathrm{vs}}_\C)$ to $W^{\mathrm{vs}}_\C$ induced by $\varphi_\C$ and $g^\circ_\C$ is the unique homorphism from $Z_\C(2\lambda_n)$ to $\varphi^{-1}_\C(W^{\mathrm{vs}}_C)$ mapping the highest weight vector $v^\circ_0$ of $Z_\C(2\lambda_n)$ onto $v_0$.

Let $K := \mathrm{ker}(f^\circ_\C) = \mathrm{ker}(g^\circ_\C)$. Then $K$ is contained in the largest proper submodule $J_\C(2\lambda_n)$ of $Z_\C(2\lambda_n)$. On the other hand $K$ has finite index in $Z_\C(2\lambda_n)$, since $\mathrm{dim}(W^{\mathrm{vs}}_\C)$ is finite. Therefore, by Weyl's theorem, if $K$ is properly contained in $J_\C(2\lambda_n)$ then $Z_\C(2\lambda_n)$ splits as a direct sum of two submodules $X$ and $Y$ where $Y = J_\C(2\lambda_n)/K$ and $X\cong Z_\C(2\lambda_n)/J_\C(2\lambda_n)\cong V_\C(2\lambda_n)$. However, if so, then $K+X$ is a proper submodule of $Z_\C(2\lambda_n)$ not contained in $J_\C(2\lambda_n)$. This contradicts the maximality of $J_\C(2\lambda_n)$. Therefore $K = J_\C(2\lambda_n)$, namely $g^\circ_\C$ induces an isomorphism from $V_\C(2\lambda_n) = Z_\C(2\lambda_n)/J_\C(2\lambda_n)$ to $\varphi_\C^{-1}(W^{\mathrm{vs}}_\C)$. In short, we may assume that $V_\C(2\lambda_n) = \varphi_\C^{-1}(W^{\mathrm{vs}}_\C)$. Thus $V_\C(2\lambda_n) \cong W^{\mathrm{vs}}_\C$, since $\varphi_\C$ is an isomorphism. This is enough to prove Theorem \ref{main theorem 1} when $\F = \C$.

\subsubsection{The general case}

Let now $\F$ be an arbitrary field. As $G_\F$ is a subgroup of $L_\F$ and the $G_\F$ orbit of $v_0$ spans $V_\F(2\lambda_n)$, we can still assume that $V_\F(2\lambda_n)$ is contained in $V_\F(2\omega_1)$. The mapping $\varphi_\F$ defined in the proof of Theorem \ref{projective veronesean} is an isomorphism from $V_\F(2\omega_1)$ to ${\cal W}_\F$. Clearly, it induces an isomorphism from the $G_\F$-submodule of $V_\F(2\omega_1)$ spanned by the $G_\F$-orbit of $v_0$ to the $G_\F$-submodule of ${\cal W}_\F$ spanned by the $G_\F$-orbit of $\bar{v}_0$. The
$G_\F$-submodule of $V_\F(2\omega_1)$ spanned by the $G_\F$-orbit of $v_0$ is the same as $V_\F(2\lambda_n)$, as noticed above, while the $G_\F$-submodule of ${\cal W}_\F$ spanned by the $G_\F$-orbit of $\bar{v}_0$ is just $W^{\mathrm{vs}}_\F$. Therefore $V_\F(2\lambda_n)\cong W^{\mathrm{vs}}_\F$. Theorem \ref{main theorem 1} is proved.

\section{Proof of Theorems \ref{main theorem 2} and \ref{main theorem 4new}}\label{sec3}

Throughout this section $\mathrm{char}(\F) = 2$. We recall that $\varepsilon^{\mathrm{vs}}_n\cong \tilde{\varepsilon}_n$, by Theorem \ref{main theorem 1}. As for notation and terminology, in general we will stick to the notation for vector spaces, but from time to time we will also use some terminology from projective geometry. The context will help avoid any ambiguity. Anyway, the symbol $\mathrm{dim (.)}$ will always denote vector dimensions, even when we are using a projective terminology. We will freely use the symbol $\Delta_n$ to denote both the point-line geometry $\Delta_n$ and its point-set. The context will make it clear in which sense that symbol is taken.

The subspaces ${\cal N}_1$, ${\cal N}_2$ and ${\cal N}_3$ of $W^{\mathrm{vs}}_n$ and the mappings $\iota_{n-1}$ and $\iota_{n-2}$ are defined as in Subsection \ref{sec1.3}. As in that subsection, for a line $l$ of $\Delta_n$ we denote by $n_l$ the nucleus of the conic $C_l := \varepsilon_n^{\mathrm{vs}}(l)$.

This section is organized as follows. We shall firstly prove Lemma \ref{main lemma n = 2} and part (1) of Theorem \ref{main theorem 4new}. Next we prove claims (1)--(4) of Theorem \ref{main theorem 2}. After that, we turn to part (2) of Theorem \ref{main theorem 4new}. We will use it to prove claim (5) of Theorem \ref{main theorem 2}.

\subsection{Proof of Lemma \ref{main lemma n = 2}}

Let $n = 2$. Then the points of $\varepsilon_2^{\mathrm{spin}}(\Delta_2)$ are the points of $\PG(3,\F)$ and the lines of $\varepsilon_2^{\mathrm{spin}}(\Delta_2)$ are the projective lines of $\PG(3,\F)$ that are totally isotropic for a non-degenerate alternating form $\alpha$ of $V(4,\F)$. We may assume that $\alpha$ admits the following expression with respect to the natural basis of $V(4,\F)$:
\[\alpha((x_1,x_2,x_3,x_4),(y_1,y_2,y_3,y_4)) ~ = ~ x_1y_4+x_2y_3+x_3y_2+x_4y_1.\]
The quadric veronesean map $\nu_4$ maps $(x_1, x_2, x_3, x_4)$ onto $(x_{i,j})_{1\leq i\leq j\leq 4}$ where $x_{i,j} = x_ix_j$. If $l$ is a line of $\PG(3,\F)$ totally isotropic for $\alpha$ then the nucleus $n_l$ of the conic $\nu_4(l)$ belongs to the $5$-dimensional subspace $S$ of $V(10,\F)$ defined by the following equations: $x_{i,i}=0 ~ \mbox{for}~ 1\leq i\leq 4 ~ \mbox{and} ~ x_{1,4}=x_{2,3}.$
It also belongs to the quadric $Q$ of $S$ defined by the following quadratic equation: $x_{1,4}^2+x_{1,2}x_{3,4}+x_{1,3}x_{2,4}=0.$
Moreover, every point of $Q$ is equal to $n_l$ for some line $l$ of $\PG(3,\F)$ totally isotropic for $\alpha$. It follows that $S = {\cal N}_2$ and $\iota_1\cong \varepsilon_1^{\mathrm{gr}}$. Lemma \ref{main lemma n = 2} is proved.

\subsection{Proof of part (1) of Theorem \ref{main theorem 4new}}\label{Lemma n = 2 e (1) Ther 4new}

Let $n = 2$ but now assume that $\F$ is perfect. We know from \cite{CP1} that $\varepsilon_2^{\mathrm{gr}}$ is a quotient of $\varepsilon_2^{\mathrm{vs}}$ ($\cong \tilde{\varepsilon}_2$ by Theorem \ref{main theorem 1}). Let $\pi$ be the projection of $\varepsilon_2^{\mathrm{vs}}$ onto $\varepsilon_2^{\mathrm{gr}}$ and $K := \mathrm{ker}(\pi)$. The space $W^{\mathrm{vs}}_2/K$ is 9-dimensional and, regarded as a $G$-module, it only contains two proper non-trivial $G$-submodules $S_1$ and $S_2$, where $S_2\subset S_1$, $\mathrm{dim}(S_1) = 5$ and $\mathrm{dim}(S_2) = 4$. Moreover, $\varepsilon_2^{\mathrm{gr}}/S_1\cong \varepsilon_2^{\mathrm{spin}}$ and $\varepsilon_2^{\mathrm{gr}}/S_2\cong \varepsilon_2^{\mathrm{sp}}$ (see \cite[Theorem 1.3]{CP1}). If ${\cal N}_3\neq K$ then $\langle {\cal N}_2, K\rangle/K$ is a 1-dimensional submodule of $W^{\mathrm{vs}}_2/K$. However $W^{\mathrm{vs}}_2/K$ does not contain any such submodule. Therefore ${\cal N}_3 = K$. It follows that $\varepsilon_2^{\mathrm{vs}}/{\cal N}_3 = \varepsilon_2^{\mathrm{gr}}$.

Since $\F$ is perfect, the embedding $\varepsilon_2^{\mathrm{sp}}$ is also a quotient of $\varepsilon_2^{\mathrm{vs}}$ (see \cite[Theorem 1.3]{CP1}). Let $\pi'$ be the projection of $\varepsilon_2^{\mathrm{vs}}$ onto $\varepsilon_2^{\mathrm{sp}}$ and $K' := \mathrm{ker}(\pi')$. Let $l$ be a line of $\Delta_2$ and put $K'_l := \langle C_l\rangle\cap K'$. Since $\varepsilon^{\mathrm{sp}}_2$ is projective, $\langle C_l\rangle/K'_l$ is 2-dimensional and the projection of $\langle C_l\rangle$ onto $\langle C_l\rangle/K'_l$ induces a bijection from $C_l$ to the set of 1-dimensional subspaces of $\langle C_l\rangle/K'_l$. This can be only if $K'_l = n_l$. Thus, ${\cal N}_2\subseteq K'$. By the previous paragraph, $K'/K = S_1$ and $K'$ is the unique 5-dimensional submodule of $W^{\mathrm{vs}}_2$. On the other hand, $\mathrm{dim}({\cal N}_2) = 5$. Hence ${\cal N}_2 = K'$, namely
$\varepsilon_2^{\mathrm{vs}}/{\cal N}_2 = \varepsilon_2^{\mathrm{sp}}$.

It is clear from the above that $W^{\mathrm{vs}}_2 \supset {\cal N}_1\supset{\cal N}_2\supset{\cal N}_3 \supset 0$ is a composition series and $\mathrm{dim}(W^{\mathrm{vs}}_2/{\cal N}_1) = \mathrm{dim}({\cal N}_2/{\cal N}_3) = 4$ and $\mathrm{dim}({\cal N}_1/{\cal N}_2) = \mathrm{dim}({\cal N}_3) = 1$. Part (1) of Theorem \ref{main theorem 4new} is proved.

\subsection{Proof of claims (1)--(4) of Theorem \ref{main theorem 2}}

Throughout this section $\F$ is assumed to be perfect.

\begin{lemma}\label{embedding properties 1}
 $\varepsilon_n^{\mathrm{vs}}(\Delta_n) \cap \mathcal{N}_1 = \emptyset$.
\end{lemma}
\pr
By way of contradiction suppose that $\varepsilon_n^{\mathrm{vs}}(X)\in \PG(\mathcal{N}_1)$ for some point $X\in\Delta_n$. By the transitivity of $G$ on $\Delta_n$ and since $G$ stabilizes ${\cal N}_1$ we obtain that
$\varepsilon_n^{\mathrm{vs}}(\Delta_n)\subseteq \PG(\mathcal{N}_1)$. Therefore, $W^{\mathrm{vs}}_n\subseteq \mathcal{N}_1.$ By \cite{TVM2004} the first $2^n$ coordinates of any vector in $\mathcal{N}_1$ are null. Hence the first $2^n$ coordinates of any vector of $W^{\mathrm{vs}}_n$ are null, but this is obviously false.  \eop

\begin{lemma}\label{embedding properties 3/2}
We have $\langle C_l\rangle \cap \mathcal{N}_1 = n_l$ for any line $l$ of $\Delta_n$.
\end{lemma}
\pr By definition, $n_l\in \mathcal{N}_1$. Suppose that $\mathrm{dim}(\langle C_l\rangle \cap \mathcal{N}_1)>1$.  Then $\langle C_l\rangle \cap \mathcal{N}_1$ contains a projective line $n$ through $n_l$. As $\F$ is perfect, $n$ is tangent to $C_l$ at a projective point $x$. Since $C_l = \varepsilon_n^{\mathrm{vs}}(l)$, we have $x = \varepsilon_n^{\mathrm{vs}}(X)$ for a point $X\in l$. However $n\subseteq{\cal N}_1$. Hence $\varepsilon_n^{\mathrm{vs}}(X)\in \varepsilon_n^{\mathrm{vs}}(\Delta_n)\cap \PG({\cal N}_1)$, a contradiction with Lemma \ref{embedding properties 1}. Therefore $\mathrm{dim}(\langle C_l\rangle \cap \mathcal{N}_1) = 1$, namely $\langle C_l\rangle \cap \mathcal{N}_1 = n_l$.  \eop

\begin{lemma}\label{lemmino banale}
We have $\langle C_l\rangle\cap\varepsilon^{\mathrm{vs}}(\Delta_n) = C_l$ for every line $l$ of $\Delta_n$.
\end{lemma}
\pr The proof is straightforward. We leave it for the reader. \eop

\begin{lemma}\label{embedding properties 2}
We have $\langle \varepsilon_n^{\mathrm{vs}}(X),\varepsilon_n^{\mathrm{vs}}(Y)\rangle \cap \mathcal{N}_1 =0$
for any two distinct points $X, Y\in \Delta_n$.
\end{lemma}
\pr By way of contradiction suppose that $\langle \varepsilon^{\mathrm{vs}}_n(X_0),\varepsilon^{\mathrm{vs}}_n(Y_0)\rangle \cap \mathcal{N}_1\not = 0$ for two distinct points $X_0, Y_0 \in \Delta_n$. By Lemma \ref{embedding properties 1} neither $\varepsilon_n^{\mathrm{vs}}(X_0)$ nor $\varepsilon_n^{\mathrm{vs}}(Y_0)$ belongs to
$\mathcal{N}_1$. Hence $z_0 := \langle \varepsilon_n^{\mathrm{vs}}(X_0),\varepsilon_n^{\mathrm{vs}}(Y_0)\rangle \cap \mathcal{N}_1$ is a projective point. Again by Lemma \ref{embedding properties 1}, $z_0\not\in \varepsilon^{\mathrm{vs}}(\Delta_n)$. Let $d = d(X_0,Y_0)$ be the distance between $X_0$ and $Y_0$ in the collinearity graph of $\Delta_n$. Suppose firstly that $d = 1$, namely $X_0$ and $Y_0$ are collinear. Let $l$ be the line of $\Delta_n$ through $X_0$ and $Y_0$. By Lemma \ref{embedding properties 3/2}, $\langle C_l\rangle\cap\mathcal{N}_1$ is the nucleus $n_l$ of the conic $C_l$. Hence $z_0 = n_l$. Consequently, the projective line $\langle \varepsilon_n^{\mathrm{vs}}(X_0),\varepsilon_n^{\mathrm{vs}}(Y_0)\rangle$ passes through the nucleus of $C_l$. Hence it is a tangent line of $C_l$. However, that line contains two points of $C_l$, namely $\varepsilon_n^{\mathrm{vs}}(X_0)$ and $\varepsilon_n^{\mathrm{vs}}(Y_0)$. We have reached a contradiction. Therefore $d > 1$.

Since the group $G$ is distance-transitive on the set of points of $\Delta_n$, the subspace $\langle \varepsilon_n^{\mathrm{vs}}(X),\varepsilon_n^{\mathrm{vs}}(Y)\rangle\cap{\cal N}_1$ is a projective point for any two points $X, Y \in \Delta_n$ at mutual distance $d$. Choose two collinear points $Y_1, Y_2\in \Delta_n$ and a point $X$ at distance $d$ from both of them and let $l$ be the line of $\Delta_n$ through $Y_1$ and $Y_2$. By Lemma \ref{lemmino banale}, $\varepsilon_n^{\mathrm{vs}}(X)$ does not belong to the projective plane spanned by the conic $C_l$ because $X\not \in l$. Hence $\mathrm{dim}(\langle\varepsilon_n^{\mathrm{vs}}(Y_1), \varepsilon_n^{\mathrm{vs}}(Y_2), \varepsilon_n^{\mathrm{vs}}(X)\rangle) = 3$. Since $d(Y_1,X) = d(Y_2,X) = d$, each of the projective lines $\langle \varepsilon_n^{\mathrm{vs}}(X), \varepsilon_n^{\mathrm{vs}}(Y_1)\rangle$ and
$\langle \varepsilon_n^{\mathrm{vs}}(X), \varepsilon_n^{\mathrm{vs}}(Y_2)\rangle$ meets ${\cal N}_1$ in a point. For $i = 1, 2$ put $z_i:=\langle \varepsilon_n^{\mathrm{vs}}(X), \varepsilon_n^{\mathrm{vs}}(Y_i)\rangle$ and let $S = \langle z_1, z_2, n_l\rangle$. Clearly $S\subseteq {\cal N}_1$.
 Hence $\PG(S)\cap C_l = \emptyset$ by Lemma \ref{embedding properties 1}. On the other hand, $S$ contains the nucleus $n_l$ of $C_l$. If $\mathrm{dim}(S\cap\langle C_l\rangle) > 1$ then $S$ contains a line tangent to $C_l$ (because $\F$ is perfect), contrary to the fact that $\PG(S)\cap C_l = \emptyset$. Hence $\mathrm{dim}(S\cap\langle C_l\rangle) = 1$, namely $S\cap\langle C_l\rangle = n_l$. However $S$ is contained in $\langle C_l, \varepsilon^{\mathrm{vs}}_n(X)\rangle$ and $\mathrm{dim}(\langle C_l, \varepsilon^{\mathrm{vs}}_n(X)\rangle) = 4$. Therefore $\mathrm{dim}(S) = 2$, namely $S$ is a projective line. Thus, $z_1$, $z_2$ and $n_l$ are collinear. It follows that the projective plane $\langle\varepsilon_n^{\mathrm{vs}}(Y_1), \varepsilon_n^{\mathrm{vs}}(Y_2), \varepsilon_n^{\mathrm{vs}}(X)\rangle$ contains $n_l$. Hence $\langle\varepsilon_n^{\mathrm{vs}}(Y_1), \varepsilon_n^{\mathrm{vs}}(Y_2), \varepsilon_n^{\mathrm{vs}}(X)\rangle = \langle C_l\rangle$, since it contains three non-collinear points of the plane $\langle C_l\rangle$, namely $\varepsilon_n^{\mathrm{vs}}(Y_1)$, $\varepsilon_n^{\mathrm{vs}}(Y_2)$ and $n_l$. Consequently $\varepsilon^{\mathrm{vs}}_n(X)\in \langle C_l\rangle$, while $\varepsilon^{\mathrm{vs}}_n(X)\not\in \langle C_l\rangle$ by Lemma \ref{lemmino banale}. We have reached a final contradiction.  \eop

\begin{prop}\label{claim (1)}
Each of ${\cal N}_1$, ${\cal N}_2$ and ${\cal N}_3$ defines a quotient of $\varepsilon^{\mathrm{vs}}_n$. Moreover:\\
(1) The embeddings $\varepsilon^{\mathrm{vs}}_n/{\cal N}_1$ and $\varepsilon_n^{\mathrm{vs}}/{\cal N}_2$ are projective while $\varepsilon_n^{\mathrm{vs}}/{\cal N}_3$ is veronesean. \\
(2) The embedding $\varepsilon_n^{\mathrm{gr}}$ is a quotient of $\varepsilon_n^{\mathrm{vs}}/{\cal N}_3$.\\
(3) For every projective embedding $\varepsilon$ of $\Delta_n$, if $\varepsilon$ is a quotient of $\varepsilon^{\mathrm{vs}}_n$ then $\varepsilon$ is also a quotient of $\varepsilon_n^{\mathrm{vs}}/{\cal N}_2$. In particular, $\varepsilon_n^{\mathrm{sp}}$ is a quotient of $\varepsilon_n^{\mathrm{vs}}/{\cal N}_2$.\\
(4) Each of the embeddings $\varepsilon^{\mathrm{vs}}_n/{\cal N}_1$, $\varepsilon_n^{\mathrm{vs}}/{\cal N}_2$ and $\varepsilon_n^{\mathrm{vs}}/{\cal N}_2$ is $G$-homogeneous.
\end{prop}
\pr By Lemmas \ref{embedding properties 1} and \ref{embedding properties 2}, the subspace ${\cal N}_1$ defines a quotient of $\varepsilon^{\mathrm{vs}}_n$. Since ${\cal N}_2\subseteq {\cal N}_1$, the subspace ${\cal N}_2$ also defines a quotient of $\varepsilon^{\mathrm{vs}}_n$. Similarly, ${\cal N}_3$ defines a quotient of $\varepsilon^{\mathrm{vs}}_n$.

By Lemma \ref{embedding properties 3/2}, for every line $l$ of $\Delta_n$ the points of the image of $l$ by $\varepsilon_n^{\mathrm{vs}}/{\cal N}_1$ bijectively correspond to the projective lines of $\langle C_l\rangle$ joining $n_l$ to points of $C_l$, namely tangents to $C_l$. The same holds for the $\varepsilon_n^{\mathrm{vs}}/{\cal N}_2$-image of $l$. However, since $\F$ is perfect, all lines of $\langle C_l\rangle$ through $n_l$ are tangents to $C_l$. Hence $\varepsilon_n^{\mathrm{vs}}/{\cal N}_1$ and $\varepsilon_n^{\mathrm{vs}}/{\cal N}_2$ map $l$ onto a line of $\PG(W^{\mathrm{vs}}_n)/{\cal N}_1$ and $\PG(W^{\mathrm{vs}}_n)/{\cal N}_2$ respectively, namely $\varepsilon^{\mathrm{vs}}_n/{\cal N}_1$ and $\varepsilon_n^{\mathrm{vs}}/{\cal N}_2$ are projective embeddings.

We know from \cite{CP1} that $\varepsilon_n^{\mathrm{gr}}$ is a quotient of $\varepsilon_n^{\mathrm{vs}}$ ($\cong \tilde{\varepsilon}_n$ by Theorem \ref{main theorem 1}). Let $\pi$ be the projection of $\varepsilon_n^{\mathrm{vs}}$ onto $\varepsilon_n^{\mathrm{gr}}$ and $K := \mathrm{ker}(\pi)$. Let $X$ be an element of $\Delta$ of type $n-2$. By part (1) of Theorem \ref{main theorem 4new} applied to the quad $Q_X$ we see that $n_{Q_X} = K\cap\langle\varepsilon_2^{\mathrm{vs}}(Q_X)\rangle$. Therefore ${\cal N}_3\subseteq K$, namely the embedding $\varepsilon_n^{\mathrm{vs}}/K \cong \varepsilon_n^{\mathrm{gr}}$ is a quotient of $\varepsilon_n^{\mathrm{vs}}/{\cal N}_3$. Since $\varepsilon_n^{\mathrm{gr}}$ is veronesean, $\varepsilon_n^{\mathrm{vs}}/{\cal N}_3$ is also veronesean.

Let $K'$ be a subspace of $W_n^{\mathrm{vs}}$ defining a quotient of $\varepsilon^{\mathrm{vs}}_n$ and assume that $\varepsilon_n^{\mathrm{vs}}/K'$ is a projective embedding. As in the second paragraph of Subsection \ref{Lemma n = 2 e (1) Ther 4new} one can see that $K'_l = n_l$ for every line $l$ of $\Delta_n$. Thus, ${\cal N}_2\subseteq K'$, namely $\varepsilon_n^{\mathrm{vs}}/K'$ is a quotient of $\varepsilon_n^{\mathrm{vs}}/{\cal N}_2$.

Since $\F$ is perfect, the embedding $\varepsilon_n^{\mathrm{sp}}$ is also a quotient of $\varepsilon_n^{\mathrm{vs}}$ (see \cite[Theorem 1.3]{CP1}). Hence $\varepsilon_n^{\mathrm{sp}}$ is a quotient of $\varepsilon_n^{\mathrm{vs}}/{\cal N}_2$, by the previous paragraph.

Finally, since each of ${\cal N}_1$, ${\cal N}_2$ and ${\cal N}_3$ is stabilized by $G$, each of $\varepsilon^{\mathrm{vs}}_n/{\cal N}_1$, $\varepsilon_n^{\mathrm{vs}}/{\cal N}_2$ and $\varepsilon_n^{\mathrm{vs}}/{\cal N}_3$ is $G$-homogeneous. \eop

\bigskip

Proposition \ref{claim (1)} covers claims (1) and (3) of Theorem \ref{main theorem 2}. We shall now turn to claim (2), but we firstly recall a few definitions on hyperplanes and projective embeddings of dual polar spaces. We recall that a {\em hyperplane} of a point-line geometry $\Gamma = ({\cal P},{\cal S})$ is a proper subset $H$ of the point-set $\cal P$ of $\Gamma$ such that for every line $l\in{\cal L}$ either $|l\cap H| = 1$ or $l\subseteq H$. Assume that $\Gamma$ is a dual polar space. Then for every point $x\in{\cal P}$ the set $H_x$ of points of $\Gamma$ at non-maximal distance from $x$ is a hyperplane of $\Gamma$. It is called a {\em singular hyperplane}.    Let $\varepsilon:\Gamma\rightarrow W$ be a projective embedding of $\Gamma$, for a vector space $W$. We say that a hyperplane $H$ of $\Gamma$ {\em arises from} $\varepsilon$ if $\varepsilon(H)$ spans a hyperplane of $W$. The embedding $\varepsilon$ is said to be {\em polarized} if every singular hyperplane of $\Gamma$ arises from $\varepsilon$.

\begin{co}\label{polarized}
Both $\varepsilon^{\mathrm{vs}}_n/{\cal N}_1$ and $\varepsilon_n^{\mathrm{vs}}/{\cal N}_2$ are polarized.
\end{co}
\pr The embeddings $\varepsilon^{\mathrm{vs}}_n/{\cal N}_1$ and $\varepsilon_n^{\mathrm{vs}}/{\cal N}_2$ are $G$-homogeneous, by Proposition \ref{claim (1)}. Hence they are polarized, by Blok, Cardinali, De Bruyn, Pasini \cite{BCDBP08}.  \eop

\begin{lemma}\label{quad embedding}
Let $n = 2$. Then $\varepsilon^{\mathrm{vs}}_2/\mathcal{N}_1 \cong \varepsilon_2^{\mathrm{spin}}$.
\end{lemma}
\pr
Let $n=2$. Then $\Delta_2$ is isomorphic to the generalized quadrangle $W(3,\F)$ associated to $\mathrm{Sp}(4,\F)$ and $\varepsilon_2^{\mathrm{spin}}$ embeds $\Delta_2$ as $W(3,\F)$ in $\PG(3,\F)$. All points of $\PG(3,\F)$ belong to $W(3,\F)$. Hence $\varepsilon^{\mathrm{vs}}_2(\Delta_2) = \nu_4(\PG(3,\F))$. Accordingly, ${\cal N}_1$ is equal to the nucleus subspace $\cal N$ of $V(10,\F)$ relative to $\nu_4$. The latter has codimension 4 in $V(10,\F)$ (Thas and Van Maldeghem \cite{TVM2004}; we warn that only finite fields of even order are considered in \cite{TVM2004}, but the above cited result, likewise most of what is said in \cite{TVM2004}, also holds for any perfect field of characteristic 2). So, $\mathrm{dim}(\varepsilon^{\mathrm{vs}}_2/{\mathcal{N}_1}) = 4$. On the other hand, the generalized quadrangle $\Delta_2$ ($\cong W(3,\F)\cong Q(4,\F)$) admits just two projective embeddings, namely $\varepsilon_2^{\mathrm{sp}}$, which is 5-dimensional, and $\varepsilon_2^{\mathrm{spin}}$, which is 4-dimensional. Thus, $\varepsilon^{\mathrm{vs}}_2/{\mathcal{N}_1} \cong \varepsilon_2^{\mathrm{spin}}$.  \eop

\bigskip

We are now ready to prove claim (2) of Theorem \ref{main theorem 2}.

\begin{prop}\label{claim (2)}
$\varepsilon_n^{\mathrm{vs}}/{\cal N}_1 \cong \varepsilon_n^{\mathrm{spin}}$.
\end{prop}
\pr We have $\mathrm{dim}(\varepsilon^{\mathrm{vs}}_n/{\cal N}_1) = 2^n$ by Corollary \ref{polarized}, Lemma \ref{quad embedding} and De Bruyn \cite[Theorem 1.6]{DB}. On the other hand, every polarized projective embedding of $\Delta_n$ has dimension at least $2^n$ by De Bruyn and Pasini \cite{DBPa2007}. Moreover, the minimal polarized projective embedding of $\Delta_n$ is unique, by Cardinali, De Bruyn and Pasini \cite{CaDBPa2007}. Since $\mathrm{dim}(\varepsilon_n^{\mathrm{vs}}/{\cal N}_1) = \mathrm{dim}(\varepsilon_n^{\mathrm{spin}}) = 2^n$, both $\varepsilon_n^{\mathrm{vs}}/{\cal N}_1$ and $\varepsilon_n^{\mathrm{spin}}$ coincide with the minimal polarized projective embedding of $\Delta_n$. Hence $\varepsilon_n^{\mathrm{vs}}/{\cal N}_1 \cong \varepsilon_n^{\mathrm{spin}}$. \eop

\bigskip

We still must prove claim (4) of Theorem \ref{main theorem 2}. We shall do it now.

\begin{prop}\label{Claim (3)}
The mapping $\iota_{n-1}$ is a projective embedding of $\Delta_{n-1}$ in ${\cal N}_2$.
\end{prop}
\pr We firstly prove that $\iota_{n-1}$ maps lines of ${\Delta}_{n-1}$ onto lines of $\PG({\cal N}_2)$. Let $l_{X,Y}:=\{Z \in \Delta_{n-1} | X\subset Z\subset Y\}$ be a line of $\Delta_{n-1}$,
 where $\{X,Y\}$ is an $\{n-2,n\}$-flag of $\Delta$. Then $l_{X,Y}$ is a line in the dual $Q^*_X$ of the quad $Q_X$, namely the $(n-1)$-grassmannian of $\mathrm{Res}^+(X)$. By Lemma \ref{main lemma n = 2}, the mapping induced by $\iota_{n-1}$ on $Q^*_X$ is a projective embedding of $Q^*_X$ in a 5-dimensional subspace of ${\cal N}_2$. Hence $\iota_{n-1}$ induces a bijection from $l_{X,Y}$ to a line of $\PG({\cal N}_2)$.

It remains to prove that $\iota_{n-1}$ is injective. By way of contradiction suppose that $\iota_{n-1}(Z_1)=\iota_{n-1}(Z_2)$ for two distinct point $Z_1$ and $Z_2$ of $\Delta_{n-1}$. Let $d = d(Z_1,Z_2)$ be the distance between $Z_1$ and $Z_2$ in the collinearity graph of $\Delta_{n-1}$. By the above, $Z_1$ and $Z_2$ are non-collinear, namely $d > 1$. Since $G$ acts distance-transitively on the set of points of $\Delta_{n-1}$, we have $\iota_{n-1}(Z)=\iota_{n-1}(Z')$ for any two points $Z, Z'\in \Delta_{n-1}$ at distance $d$. Take now two collinear points $Z_1$ and $Z_2$ of $\Delta_{n-1}$. Then $\iota_{n-1}(Z_1) \neq \iota_{n-1}(Z_2)$, by the above. On the other hand, there exists a point $Z_3\in\Delta_{n-1}$ such that $d(Z_3,Z_1) =
d(Z_3,Z_2) = d$. Therefore $\iota_{n-1}(Z_1) = \iota_{n-1}(Z_3) = \iota_{n-1}(Z_2)$, whence $\iota_{n-1}(Z_1) = \iota_{n-1}(Z_2)$. We have reached a final contradiction.  \eop

\subsection{Proof of part (2) of Theorem \ref{main theorem 4new}}

Let $n = 3$ and $\F$ be perfect. We recall that, in view of Theorem \ref{main theorem 1}, $\varepsilon^{\mathrm{vs}}_3\cong \tilde{\varepsilon}_3$, namely $W^{\mathrm{vs}}_3\cong V(2\lambda_3)$. The isomorphisms $\iota_2\cong \tilde{\varepsilon}_2$, $\iota_1\cong\tilde{\varepsilon}_1$, $\tilde{\varepsilon}_3/{\cal N}_2\cong \varepsilon^{\mathrm{sp}}_3$ and $\tilde{\varepsilon}_3/{\cal N}_3\cong \varepsilon^{\mathrm{gr}}_3$ can be proved by straightforward calculations, writing down explicit expressions for $\varepsilon^{\mathrm{sp}}_3$ and $\varepsilon^{\mathrm{vs}}_3$ and computing coordinates or equations for all objects involved here. However, these calculations are too cumbersome to be exposed in a perspicuous way. We prefer to adopt a different strategy, exploiting more syntetic arguments as far as possible.

According to \cite[Theorem 1.3]{CP1} and claim (2) of Theorem \ref{main theorem 2} the module $V(2\lambda_3)$ ($\cong W^{\mathrm{vs}}_3$ by Theorem \ref{main theorem 1}) admits a series of submodules $V(2\lambda_3) \supset {\cal N}_1 \supset K_2 \supset K_1 \supset 0$ where $\mathrm{dim}(V(2\lambda_3)/{\cal N}_1) = 8$, $\mathrm{dim}({\cal N}_1/K_2) = 6$, $\mathrm{dim}(K_2/K_1) = 14$ and $\mathrm{dim}(K_1) = 7$
and such that $\tilde{\varepsilon}_3/{\cal N}_1 \cong \varepsilon_3^{\mathrm{spin}}$, $\tilde{\varepsilon}_3/K_2\cong \varepsilon^{\mathrm{sp}}_3$ and $\tilde{\varepsilon}_3/K_1 \cong \varepsilon_3^{\mathrm{gr}}$. Moreover, $K_2/K_1$ hosts $\varepsilon^{\mathrm{sp}}_2$. Note also that ${\cal N}_1/K_2$ hosts $\varepsilon_1^{\mathrm{sp}}$ (see Cardinali and Lunardon \cite{CL}).

\begin{lemma}\label{K2/K1}
The section $K_2/K_1$ is an irreducible $G$-module.
\end{lemma}
\pr As noticed above, $K_2/K_1$ hosts $\varepsilon_2^{\mathrm{sp}}$. Namely $K_2/K_1$ is isomorphic to the 14-dimensional module $\overline{V}_2$ for $\mathrm{Sp}(6,\F)$ ($\cong G$) contained in the exterior square $\overline{V}\wedge\overline{V}$ of $\overline{V} = V(6,\F)$ (see the proof of Theorem 1.3 of \cite{CP1}). The module $\overline{V}_2$ is irreducible, as one can see by using the main result of Premet and Suprunenko \cite{PS}, which also holds in characteristic 2, by Baranov and Suprunenko \cite{BS}.  \eop

\begin{lemma}\label{dim(N2) = 21}
We have ${\cal N}_2 = K_2 \cong V(\lambda_2)$.
\end{lemma}
\pr We have ${\cal N}_2\subseteq K_2$ by claim (3) of Theorem \ref{main theorem 2}. On the other hand, ${\cal N}_2\not\subseteq K_1$, since $\tilde{\varepsilon}_3/{\cal N}_2$ is a projective embedding (claim (3) of Theorem \ref{main theorem 2}) while $\tilde{\varepsilon}_3/K_1 \cong \varepsilon_3^{\mathrm{gr}}$, which is veronesean. Therefore, if ${\cal N}_2$ contains $K_1$ then ${\cal N}_2 = K_2$ by Lemma \ref{K2/K1}.

In order to go on we need to recover $K_1$ inside $V(2\lambda_3)$. Let $\mathfrak{L}$ be the Lie algebra of $G = \mathrm{Spin}(7,\F)$ and let $\alpha_1, \alpha_2, \alpha_3$ be the three simple roots of the root system of type $B_3$. We recall that $2\lambda_3 = \alpha_1 + 2\alpha_2 + 3\alpha_3$ (see e.g. Humphreys \cite{Humphreys}). Let $v_0$ be a vector of $V(2\lambda_3)$ of weight $2\lambda_3$ and let $x_{-\alpha_3}$ be a non-zero element of $X_{-\alpha_3}$, where $X_{-\alpha_3}$ is the 1-dimensional subspace of $\mathfrak{L}$ associated to $-\alpha_3$. Let $v_1 := x_{-\alpha_3}v_0$. Then $v_1$ has weight $\lambda_2 = \alpha_1 + 2\alpha_2 + 2\alpha_3$. It is straightforward to check that the submodule $\langle G(v_1)\rangle$ of $V(2\lambda_3)$ spanned by the $G$-orbit $G(v_1)$ of $v_1$ is isomorphic to $V(\lambda_2)$. Whence it is 21-dimensional. One can also check that $V(2\lambda_3)$ does not contain any 21-dimensional $G$-submodule different from $\langle G(v_1)\rangle$. On the other hand, $K_2$ is indeed a 21-dimensional $G$-submodule of $V(2\lambda_3)$. Hence $K_2 = \langle G(v_1)\rangle\cong V(\lambda_2)$.

Regarded $K_2$ as a copy of $V(\lambda_2)$, the vector $v_1$ is a vector of $V(\lambda_2)$ of weight $\lambda_2 = \alpha_2 + 2\alpha_2+ 2\alpha_3$. Let $v_2 := x_{-\alpha_3}x_{-\alpha_2}v_1$. Then $v_2$ has weight $\lambda_1 = \alpha_1 + \alpha_2 + \alpha_3$. It is straightforward to check that the submodule $\langle G(v_2)\rangle$ of $V(\lambda_2)$ spanned by the $G$-orbit $G(v_2)$ of $v_2$ is isomorphic to $V(\lambda_1)$. Whence it is 7-dimensional. One can also check that $V(\lambda_2)$ does not contain any 7-dimensional $G$-submodule different from $\langle G(v_2)\rangle$. On the other hand, $K_1$ is indeed a 7-dimensional $G$-submodule of $K_2 \cong V(\lambda_2)$. Hence $K_1 = \langle G(v_2)\rangle\cong V(\lambda_1)$. Finally, $K_1 \cong V(\lambda_1)$ admits a unique 1-dimensional submodule $K_0$. We have previously remarked that $K_2/K_1$ is irreducible. Clearly, $K_1/K_0$ is irreducible as well. Note that $K_2/K_1$ hosts $\varepsilon_2^{\mathrm{sp}}$ (see \cite[Theorem 1.3]{CP1}) while $K_1/K_0$ hosts $\varepsilon_1^{\mathrm{sp}}$.

Suppose that ${\cal N}_2\neq K_2$. By the first paragraph of this proof, ${\cal N}_2\cap K_1$ is different from either $K_1$ and ${\cal N}_2$. Moreover ${\cal N}_2+K_1 = K_2$ since $K_2/K_1$ is irreducible. As $K_1/K_0$ is irreducible, ${\cal N}_2\cap K_1\subseteq K_0$.

By (4) of Theorem \ref{main theorem 2} the mapping $\iota_2$ is a projective embedding of $\Delta_2$
in ${\cal N}_2$. By \cite[Theorem 1.5(1)]{CP1}, the embedding $\widetilde{\varepsilon}_2:\Delta_2\rightarrow V(\lambda_2)$ is universal. Therefore $\iota_2$ is a quotient of $\tilde{\varepsilon}_2$, namely ${\cal N}_2$ is isomorphic to a quotient of $V(\lambda_2)$. Thus $K_2$ admits a submodule $N$ such that $K_2/N\cong {\cal N}_2$.

Suppose that
${\cal N}_2\cap K_0 = K_0$. Then $N$ is contained in a further submodule $N' \supset N$ such that $N$ has codimension 1
in $N'$. Clearly, $N'/N$ corresponds to $K_0\subset {\cal N}_2$. On the other hand, since ${\cal N}_2+K_1 = K_2$, we have that $\mathrm{dim}({\cal N}_2) = 15$ ($= 14 + 1 = \mathrm{dim}(K_2/K_1) + \mathrm{dim}(K_0)$). Hence $\mathrm{dim}(N) = 6$. As $K_2/K_1$ is irreducible, $N$ can only sit in $K_1$. However $K_1 = V(\lambda_1)$ does not admit any 6-dimensional submodule. We have reached a contradiction, which forces us to conclude that
${\cal N}_2\cap K_0 = 0$. Thus, $K_2 \cong V(\lambda_2)$ splits as a direct sum ${\cal N}_2\oplus K_1$, where $K_1 = V(\lambda_1)$. However this is false, as one can check. We have reached a final contradiction. Therefore ${\cal N}_2 = K_2$   \eop

\begin{prop}\label{Th 4 (2) - 1}
We have $\iota_2\cong\tilde{\varepsilon}_2$ and $\varepsilon_3^{\mathrm{vs}}/{\cal N}_2\cong \varepsilon_3^{\mathrm{sp}}$.
\end{prop}
\pr Both claim immediately follows from Lemma \ref{dim(N2) = 21}.  \eop

\begin{lemma}\label{dim(N3) = 7}
${\cal N}_3 = K_1 \cong V(\lambda_1)$.
\end{lemma}
\pr We keep the notation introduced in the proof of Lemma \ref{dim(N2) = 21}. We have ${\cal N}_3\subseteq K_1$ by (3) of Theorem \ref{main theorem 2}. As $K_1/K_0$ is irreducible and does not split as a direct sum of two submodules, either ${\cal N}_3 = K_1$ or ${\cal N}_3 = K_0$.

We know that $\iota_2 \cong \tilde{\varepsilon}_2$ (Proposition \ref{Th 4 (2) - 1}). Accordingly, for a 1-element $X$ of $\Delta$ we can recover $n_{Q_X} = \iota_1(X)$ in $K_1 = V(\lambda_1)$ as the nucleus of the quadric $\tilde{\varepsilon}_2(Q^*_X)$, where $Q^*_X$ is the dual of $Q_X$, namely $Q^*_X = \mathrm{Res}(X)$ with elements of type 2 and 3 taken as points and lines respectively. It is implicit in the proof of Lemma \ref{dim(N2) = 21} that $\varepsilon_2^{\mathrm{gr}} \cong \tilde{\varepsilon}_2/K_0$. The embeddings $\tilde{\varepsilon}_2$ and $\tilde{\varepsilon}_2/K_0$ induce the same embedding on $Q^*_X$. Hence
$\varepsilon_2^{\mathrm{gr}}(Q^*_X)$ is still a quadric. Its nucleus is equal to $\langle n_{Q_X}, K_0\rangle/K_0$. It follows that $n_{Q_X}\neq K_0$. Therefore ${\cal N}_3 = K_1$.   \eop

\begin{prop}\label{Th 4 (2) - 2}
We have $\iota_1\cong\tilde{\varepsilon}_1$ and $\varepsilon_3^{\mathrm{vs}}/{\cal N}_3\cong \varepsilon_3^{\mathrm{gr}}$.
\end{prop}
\pr The isomorphism $\varepsilon_3^{\mathrm{vs}}/{\cal N}_3\cong \varepsilon_3^{\mathrm{gr}}$ immediately follows from the fact that ${\cal N}_3 = K_1$ (Lemma \ref{dim(N3) = 7}). Turning to $\iota_1$, recall that the points of $\Delta_1$ bijectively correspond to the 1-dimensional subspaces of $V(\lambda_1)$ contained in a $G$-orbit of a given highest weight vector $v_2$ of $V(\lambda_1)$ (the symbol $v_2$ is inherited from the proof of Lemma \ref{dim(N2) = 21}). Let $A$ be the point of $\Delta_1$ corresponding to $\langle v_2\rangle$, put $Q^*_A  := \mathrm{Res}(A)$ (notation as in the proof of Lemma \ref{dim(N3) = 7}) and let $n_A$ be the nucleus of the quadric $\tilde{\varepsilon}_2(Q^*_A)$. Let $G_A$, $G_{\langle v_1\rangle}$ and $G_{n_A}$ be the stabilizers of
$A$, $\langle v_1\rangle$ and $n_A$ respectively in $G$. We have $G_A = G_{\langle v_2\rangle}$ by our choice of $A$. On the other hand $G_A\leq G_{n_A}$, whence $G_A = G_{n_A}$ by the maximality of $G_A$ in $G$. Thus $G_A = G_{n_A} = G_{\langle v_2\rangle}$. Accordingly, $G_{g(A)} = G_{n_{g(A)}} = G_{\langle g(v_2)\rangle}$ for any $g\in G$. So, the points of $\PG(V(\lambda_1))$ in the $G$-orbit of $\langle v_2\rangle$ bijectively correspond to the points of $\Delta_1$ as well as to those of $\iota_1(\Delta_1)$ in such a way that a point $\langle g(v_2)\rangle$ of $\PG(V(\lambda_1))$ with $g\in G$ and a point of $\Delta_1$ or $\iota_1(\Delta_1)$ correspond if and only if they are stabilized by the same subgroup of $G$. However, $\iota_1(\Delta_1)$ sits in $\PG(K_1)$ and the $G$-modules $K_1$ and $V(\lambda_1)$ are isomorphic. It follows that, for every line $l$ of $\Delta_1$, the image $\iota_1(l)$ of $l$ is mapped by the isomorphism $K_1\cong V(\lambda_1)$ onto a line of $\PG(V(\lambda_1))$. Since isomorphism of modules are realized by means of invertible linear transformations, $\iota_1(l)$ must be a line in $\PG(K_1)$ too. Thus, $\iota_1\cong \tilde{\varepsilon}_1$.   \eop

Along the way, we have also proved that $V(2\lambda_3)\supset {\cal N}_1\supset {\cal N}_2\supset {\cal N}_3\supset K_0\supset 0$ is a composition series for the $G$-module $V(2\lambda_3)$. The proof of (2) of Theorem~\ref{main theorem 4new} is complete.

\subsection{Proof of claim (5) of Theorem \ref{main theorem 2}}

Let $n \geq 3$ and $\F$ be perfect. When $n = 3$ claim (5) of Theorem \ref{main theorem 2} is contained in part (2) of Theorem \ref{main theorem 4new}. Assume that $n > 3$. Given an $\{n-3,n-1\}$-flag $\{X,Y\}$ of $\Delta$,
 let $l_{X,Y}$ be the line of $\Delta_{n-2}$ determined by $\{X,Y\}$, namely $l_{X,Y} = \{Z\in \Delta_{n-2} | X\subset Z \subset Y\}$. The upper residue $\mathrm{Res}^+(X)$ of $X$ is a building of type $B_3$, with $\{n-2, n-1, n\}$ as the set of types and $l_{X,Y}$ is a line of its $(n-2)$-grassmannian. So, we can apply part (2) or Theorem \ref{main theorem 4new} to $\mathrm{Res}^+(X)$, obtaining that the mapping induced by $\iota_{n-2}$ on $\mathrm{Res}^+(X)$ maps $l_{X,Y}$ bijectively onto a projective line of $\PG({\cal N}_3)$.

In order to prove that $\iota_{n-2}$ is a projective embedding of $\Delta_{n-2}$ in ${\cal N}_3$ it remains to show that $\iota_{n-2}$ is injective. This can be proved by the same argument used in the proof of Proposition \ref{Claim (3)} to show that $\iota_{n-1}$ is injective. We leave the details for the reader.

\section{Universality} \label{sec4}

\subsection{Veronesean embeddings and 3-generation}\label{3-gen}

Let $\Gamma = ({\cal P},{\cal L})$ be a point-line geometry. We say that a subset $S$ of the point-set $\cal P$ of $\Gamma$ is a $3$-{\em subspace} of $\Gamma$ if $S$ contains every line $l$ of $\Gamma$ such that $|l\cap S| \geq 3$. Clearly, intersections of 3-subspaces are 3-subspaces. So, we can consider the $3$-{\em span} $\langle X\rangle^{(3)}_\Gamma$ of a subset $X\subseteq {\cal P}$, defined as the smallest 3-subspace containing $X$. We say that $X$ $3$-generates $\Gamma$ if $\langle X\rangle^{(3)}_\Gamma = {\cal P}$. The $3$-{\em generating rank} $\mathrm{grk}_3(\Gamma)$ of $\Gamma$ is the size of a smallest 3-generating set of $\Gamma$. Clearly, if $\Gamma$ admits a veronesean embedding $\nu$ then $\mathrm{dim}(\nu) \leq \mathrm{grk}_3(\Gamma)$.

\begin{prop}\label{lines of size 3}
Suppose that all lines of $\Gamma$ have exactly three points. Then $\mathrm{grk}_3(\Gamma) = |{\cal P}|$.
\end{prop}
\pr Clearly $\mathrm{grk}_3(\Gamma) \leq |{\cal P}|$. On the other hand, let $V = V(|{\cal P}|,2)$, let $B$ be a basis of $V$ and $\nu$ a bijection from $\cal P$ to $B$. Then $\nu$ is a veronesean embedding of $\Gamma$. Hence $\mathrm{grk}_3(\Gamma) = |{\cal P}|$.  \eop

\bigskip

Note that the embedding $\nu$ mentioned in the proof of Proposition \ref{lines of size 3} is absolutely universal (in the class of all veronesean or projective embeddings of $\Gamma$).

In the next proposition we consider a situation completely different from that of Proposition \ref{lines of size 3}.

\begin{prop}\label{quasi-ver-proj}
Let $\Gamma = \PG(d,\F)$ with $\F\neq \F_2$. Then $\mathrm{grk}_3(\Gamma) = {{d+2}\choose 2}$.
\end{prop}
\pr The claim that $\mathrm{grk}_3(\PG(d,\F)) \leq {{d+2}\choose 2}$ is essentially the same as \cite[Lemma 4.11]{CP1} while the inequality $\mathrm{grk}_3(\PG(d,\F)) \geq {{d+2}\choose 2}$ follows from the fact that the usual quadric veronesean map of $\PG(d,\F)$ has dimension equal to ${{d+2}\choose 2}$. These two inequalities combined together yield the equality $\mathrm{grk}_3(\PG(d,\F)) = {{d+2}\choose 2}$.  \eop

\bigskip

The next result has been proved by Thas and Van Maldeghem \cite[Theorem 1.6]{TVM2004} under the assumption that $\F$ is finite, but it is also an immediate consequence of Proposition \ref{quasi-ver-proj}. Moreover, it holds for any field $\F \neq \F_2$.

\begin{co}\label{veronesean universal}
Let $\F\neq \F_2$. Then the quadric veronesean map of $\PG(d,\F)$ is relatively universal.
\end{co}
{\bf Remarks.} (1) In Thas and Van Maldeghem \cite{TVM} embeddings $\varepsilon:\mathrm{PG}(d,\F)\rightarrow V$ are introduced, called {\em generalized veronesean embeddings}, where for every line $l$ the image $\varepsilon(l)$ of $l$ under $\varepsilon$ spans a projective plane of $\mathrm{PG}(V)$. Clearly, veronesean embeddings as defined in this paper are special cases of generalized veronesean embeddings. A complete classification of generalized veronesean embeddings is achieved in \cite{TVM} under the additional assumption that $\mathrm{dim}(V) = m \geq {{d+2}\choose 2}$ and a quite natural additional assumption when $\F$ is infinite. As a by-product of that classification, $m = {{d+2}\choose 2}$. This particular result also follows from our Proposition \ref{quasi-ver-proj}.\\
(2) Let $\F \neq \F_2$. In view of Thas and Van Maldeghem \cite{TVM} the quadric veronesean map of $\PG(d,\F)$ is absolutely universal (in the class of all veronesean embeddings of $\PG(d,\F)$) if and only if every veronesean embedding of $\PG(d,\F)$ of dimension less than ${{d+2}\choose 2}$, if any, is a quotient of a ${{d+2}\choose 2}$-dimensional veronesean embedding.

\subsection{Proof of Theorem \ref{main theorem univ1}}

Let $n = 2$. Then $\Delta_2$ is isomorphic to the generalized quadrangle $W(3,\F)$ of symplectic type and $\varepsilon^{\mathrm{spin}}_2$ is just the natural embedding of $W(3,\F)$ in $\PG(3,\F)$. Thus, $\varepsilon^{\mathrm{vs}}_2$ is nothing but the veronesean embedding of $W(3,\F)$ induced by the quadratic veronesean map $\nu_4:V(4,\F)\rightarrow V(10,\F)$, where $W(3,\F)$ is regarded as a subgeometry of $\PG(3,\F)$. As for the group $G = \mathrm{Spin}(5,\F)$, we may regard $G$ as the same group as $\mathrm{Sp}(4,\F)$, since $\mathrm{Spin}(5,\F)\cong \mathrm{Sp}(4,\F)$.

Let $\F = \F_q$ for a prime power $q > 2$. We assume $q > 2$ in view of Proposition \ref{lines of size 3}. Throughout this subsection we use the following shortened notation. We write $\Gamma$ for $W(3,q)$, $\nu$ for $\nu_4$ and we denote by $\varepsilon$ the veronesean embedding of $\Gamma$ induced by $\nu$, namely $\varepsilon = \nu\circ\iota$ where $\iota$ is the inclusion of $\Gamma$ in $\PG(3,q)$. Let $\perp$ be the collinearity relation of $\Gamma$. We recall that, given a non-collinear pair of points $\{a,b\}$ of $\Gamma$ the set $\{a,b\}^{\perp\perp}$ is called a {\em hyperbolic line} of $\Gamma$. In order to avoid any risk of confusion, we shall call the lines of $\Gamma$ {\em isotropic lines}.
We recall that $\mathrm{PG}(3,q)$ and $\Gamma$ have the same set of points. The lines of $\mathrm{PG}(3,q)$ are the isotropic and hyperbolic lines of $\Gamma$.

The embedding $\varepsilon$ is a $G$-homogeneous embedding. Let $\tilde{\varepsilon}:\Gamma\rightarrow\widetilde{V}$ be another $G$-homogeneous embedding of $\Gamma$ such that a linear mapping $\pi:\widetilde{V}\rightarrow V(10,q)$ exists satisfying the following: $\pi\circ \tilde{\varepsilon} \cong\varepsilon$ and $G$ stabilizes the kernel $K := \mathrm{ker}(\pi)$ of $\pi$. Clearly, the linear hull of $\varepsilon$ satisfies the above conditions.

In order to prove Theorem \ref{main theorem univ1} we only must prove that if $q$ is odd and greater than 3 then $\pi$ is an isomorphism, namely $K = 0$. In other words, if $K \neq 0$ then either $q$ is even or $q = 3$. Thus, for the rest of this section we assume that $K \neq 0$.

We have $\mathrm{dim}(\varepsilon) = 10$ since $\Gamma$ and $\PG(3,q)$ have the same set of points and $\mathrm{dim}(\nu) = 10$. Hence $\mathrm{dim}(\widetilde{V}) = \mathrm{dim}(K)+10 >
10$ by the assumption that $K \neq 0$.

For every isotropic line $l$ of $\Gamma$ the projection $\pi$ induces an isomorphism from $\langle\tilde{\varepsilon}(l)\rangle$ to $\langle\varepsilon(l)\rangle$. Hence $\tilde{\varepsilon}(l)$ is a conic for every isotropic line $l$ of $\Gamma$. For a hyperbolic line $l$ of $\Gamma$, put $d := \mathrm{dim}(\langle\tilde{\varepsilon}(l)\rangle)$. As $\tilde{\varepsilon}$ is $G$-homogeneous by assumption and $G$ is transitive on the set of hyperbolic lines of $\Gamma$, the dimension $d$ does not depend on the choice of the hyperbolic line $l$. Clearly, $d \leq q+1$. Moreover $d \geq 3$, since the projection of $\tilde{\varepsilon}$ onto $\varepsilon$ maps $\langle\tilde{\varepsilon}(l)\rangle$ onto $\langle\varepsilon(l)\rangle$ and $\varepsilon(l)$ is a conic in $\mathrm{PG}(9,q)$.

\begin{lemma}\label{XXX.1}
If $d = 3$ then $\tilde{\varepsilon} = \varepsilon$.
\end{lemma}
\pr If $d = 3$ then $\tilde{\varepsilon}$ is also a veronesean embedding of $\PG(3,q)$. Hence
$\tilde{\varepsilon} = \varepsilon$ by Corollary \ref{veronesean universal}.  \eop

\bigskip

By Lemma \ref{XXX.1}, if $d = 3$ then Theorem \ref{main theorem univ1} is proved. Thus, henceforth we assume $d > 3$.

The group $G$ has two orbits on the set of pairs $\{l,m\}$ of hyperbolic lines, depending on whether $l\cap m$ is either empty or a point. Hence for a pair $\{l,m\}$ of hyperbolic lines
the dimension of $\langle \tilde{\varepsilon}(l)\rangle\cap\langle\tilde{\varepsilon}(m)\rangle$
only depends on whether $l\cap m$ is empty or a point. In the sequel we are mainly interested in the case where $l\cap m$ is a point. For a pair $\{l,m\}$ of hyperbolic lines intersecting in a point we put $\delta := \mathrm{dim}(\langle \tilde{\varepsilon}(l)\rangle\cap\langle\tilde{\varepsilon}(m)\rangle) - 1.$
\begin{lemma}\label{XXX.2}
Let $a$ be a point of $\Gamma$ and $d' := \mathrm{dim}(\langle\tilde{\varepsilon}(a^\perp)\rangle)$. Then one of the following holds. \\
(1) $d' = 2d-\delta$.\\
(2) $d = q+1$ and $d' = 2d-\delta + 1$.
\end{lemma}
\pr
Chosen two distinct hyperbolic lines $l$ and $m$ in $a^\perp$, the point $a$ does not belong to $l\cup m$ while $l\cap m$ is a point, say $c$. Let $a_1,..., a_{d-1}$ and $b_1,..., b_{d-1}$ be points of $l$ and $m$ respectively such that $\{\tilde{\varepsilon}(a_1),...,\tilde{\varepsilon}(a_{d-1}), \tilde{\varepsilon}(c)\}$ is a basis of $\langle\tilde{\varepsilon}(l)\rangle$ and
$\{\tilde{\varepsilon}(b_1),...,\tilde{\varepsilon}(b_{d-1}), \tilde{\varepsilon}(c)\}$ is a basis of $\langle\tilde{\varepsilon}(m)\rangle$. Then the subspace
\[X_{a,l,m} = \langle \tilde{\varepsilon}(a), \tilde{\varepsilon}(a_1),...,\tilde{\varepsilon}(a_{d-1}), \tilde{\varepsilon}(b_1),...,\tilde{\varepsilon}(b_{d-1}), \tilde{\varepsilon}(c)\rangle\]
contains $\tilde{\varepsilon}(l)\cup\tilde{\varepsilon}(m)\cup\{\tilde{\varepsilon}(a)\}$. Every hyperbolic line of $\Gamma$ is equal to $\{a,b\}^{\perp\perp}$ for a pair $\{a,b\}$ of non-collinear points of $\Gamma.$  The conic $\tilde{\varepsilon}(n)$ is fully contained in $X_{a,l,m}$ for every isotropic line $n$ through $a$ different from the line $ac$ through $a$ and $c.$ Indeed, the line $n$ meets $l\cup m\cup \{a\}$ in three distinct points hence $\tilde{\varepsilon}(n)$ is a conic contained in $X_{a,l,m}.$ If follows that $X_{a,l,m}$ contains the images of all points of $a^\perp$ but possibly those that belong to
the line $ac$ through $a$ and $c$ and are different from either $a$ or $c$. Let $[a^\perp]_c = a^\perp\setminus(ac\setminus\{a,c\})$. Every hyperbolic line contained in $a^\perp$ but not containing $c$ meets $[a^\perp]_c$ in $q$ points. Therefore if $d < q+1$ every such hyperbolic line is mapped by $\tilde{\varepsilon}$ onto a subset of $X_{a,l,m}$. In this case $X_{a,l,m}$ also contains $\tilde{\varepsilon}(ac)$. It follows that $2d-\delta-1 \leq d' \leq 2d-\delta$. Suppose that $d' = 2d-\delta-1$. Then $\tilde{\varepsilon}(a) \in \langle \tilde{\varepsilon}(l)\cup\tilde{\varepsilon}(m)\rangle$. By projecting $\tilde{\varepsilon}$ onto $\varepsilon$ we get $\varepsilon(a) \in \langle \varepsilon(l)\cup\varepsilon(m)\rangle$. However, this is false, as one can see by a direct computation. Therefore $d' = 2d-\delta$, as in case (1) of the lemma.

Let $d = q+1$. Then we need only one point of $ac\setminus\{a,c\}$ besides $a$ and $c$ in order to span $\langle\tilde{\varepsilon}(ac)\rangle$. If that point is not mapped by $\tilde{\varepsilon}$ into $X_{a,l,m}$ then we are in case (2). \eop

\bigskip

Let $a$ and $b$ be two non-collinear points of $\Gamma$. The dimension $\mathrm{dim}(\langle \tilde{\varepsilon}(a^\perp)\rangle\cap\langle \tilde{\varepsilon}(b^\perp)\rangle)$ does not depend on the choice of the non-collinear pair $\{a,b\}$, since $G$ is transitive on the set of pairs of non-collinear points of $\Gamma$. Note that $a^\perp\cap b^\perp$ is a hyperbolic line. Hence $\mathrm{dim}(\langle \tilde{\varepsilon}(a^\perp)\rangle\cap\langle \tilde{\varepsilon}(b^\perp)\rangle ) \geq d$.

We put $\delta' :=
\mathrm{dim}(\langle \tilde{\varepsilon}(a^\perp)\rangle\cap\langle\tilde{\varepsilon}(b^\perp)\rangle )- d$. Then,
$\mathrm{dim}(\langle\tilde{\varepsilon}(a^\perp\cup b^\perp)\rangle) = 2d'-d-\delta'$. By Lemma
\ref{XXX.2}, either $2d'-d-\delta' = 3d-2\delta-\delta'$ or $2d'-d-\delta' = 3d-2\delta-\delta'+2$, the latter case occurring only if $d = q+1$.

\begin{lemma}\label{XXX.3}
One of the following holds: \\
(1) $d' = 2d-\delta$ and $\mathrm{dim}(\widetilde{V}) = 3d-2\delta-\delta' + 1$;\\
(2) $d = q+1$, $d' = 2d-\delta+1$ and $3d-2\delta-\delta' + 3\leq \mathrm{dim}(\widetilde{V}) \leq 3d-2\delta-\delta' + 4$.
\end{lemma}
\pr
Let $a$, $b$ and $c$ be three pairwise non-collinear points of $\Gamma$, with $c\not\in \{a,b\}^{\perp\perp}$. Put $X_{a,b,c} = \langle \tilde{\varepsilon}(a^\perp\cup b^\perp\cup\{c\})\rangle$. Then $2d'-d-\delta'\leq \mathrm{dim}(X_{a,b,c})\leq 2d'-d-\delta'+1$. If $\mathrm{dim}(X_{a,b,c}) = 2d'-d-\delta'$ then $\tilde{\varepsilon}(c)\in \langle\tilde{\varepsilon}(a^\perp\cup b^\perp)\rangle$. By projecting $\tilde{\varepsilon}$ onto $\varepsilon$ we get $\varepsilon(c)\in \langle\varepsilon(a^\perp\cup b^\perp)\rangle$. However a direct computation shows that this is not the case. Hence $\mathrm{dim}(X_{a,b,c}) = 2d'-d-\delta'+1$.

As $c\not\in\{a,b\}^\perp\cup\{a,b\}^{\perp\perp}$, $c^\perp$ meets $\{a,b\}^\perp$ in a single point, say $c'$. All isotropic lines through $c$ different from the line $cc'$ through $c$ and $c'$ meet $a^\perp\cup b^\perp$ in distinct points. Hence all isotropic lines through $c$ different from $cc'$ meet $a^{\perp}\cup b^{\perp} \cup \{c\}$ in three distinct points. Hence $\tilde{\varepsilon}$ maps them onto conics contained in $X_{a,b,c}$. Suppose firstly that $d' = 2d-\delta$. Then, arguing as in the proof of Lemma \ref{XXX.2} one can prove that $\tilde{\varepsilon}(c^\perp)\subseteq X_{a,b,c}$. Let now $l$ be an isotropic line meeting $a^\perp$, $b^\perp$ and $c^\perp$ in pairwise distinct points. Then the conic $\tilde{\varepsilon}(l)$ meets $X_{a,b,c}$ in three distinct points. Hence $\tilde{\varepsilon}(l)\subseteq X_{a,b,c}$. Every point $x\not\in \{a,b\}^{\perp\perp}\cup\{b,c\}^{\perp\perp}\cup\{c,a\}^{\perp\perp}$ belongs to a line $l$ as above. Hence every such point is mapped by $\tilde{\varepsilon}$ onto a point of $X_{a,b,c}$. Thus $X_{a,b,c}$ contains the images of all points of $\Gamma$ except possibly those of $(\{a,b\}^{\perp\perp}\cup\{b,c\}^{\perp\perp}\cup\{c,a\}^{\perp\perp})\setminus(a^{\perp}\cup b^{\perp}\cup c^{\perp})$. Let $x\in (\{a,b\}^{\perp\perp}\cup\{b,c\}^{\perp\perp}\cup\{c,a\}^{\perp\perp})\setminus(a^{\perp}\cup b^{\perp}\cup c^{\perp})$. As $q > 2$ by assumption, every isotropic line through $x$ contains at least three points that do not belong to $(\{a,b\}^{\perp\perp}\cup\{b,c\}^{\perp\perp}\cup\{c,a\}^{\perp\perp})$. By the above, that line is mapped by $\tilde{\varepsilon}$ onto a conic of $X_{a,b,c}$. Hence $\tilde{\varepsilon}(x)\in X_{a,b,c}$. It follows that $X_{a,b,c} = \widetilde{V}$.

Finally, let $d' = 2d-\delta+1$ (whence $d = q+1$). Then $\mathrm{dim}(X_{a,b,c}) = 3d-2\delta+3-\delta'$. The subspace $X_{a,b,c}$ contains the images of all points of $c^\perp$, except possibly those of $cc'\setminus\{c,c'\}$. If $X_{a,b,c}\supseteq \tilde{\varepsilon}(c^\perp)$ then we obtain that $\mathrm{dim}(\widetilde{V}) = 3d-2\delta+3-\delta'$. Otherwise, let $c''\in cc'\setminus\{c,c'\}$. Then $\tilde{\varepsilon}(cc')\subseteq \langle \tilde{\varepsilon}(c), \tilde{\varepsilon}(c'), \tilde{\varepsilon}(c'')$. Therefore $\langle X_{a,b,c}\cup\{\tilde{\varepsilon}(c'')\}\rangle$ contains $\tilde{\varepsilon}(c^\perp)$. As in the previous paragraph, one can prove that $\langle X_{a,b,c}\cup\{\tilde{\varepsilon}(c'')\}\rangle$ contains all of $\tilde{\varepsilon}(\Gamma)$. Hence $\langle X_{a,b,c}\cup\{\tilde{\varepsilon}(c'')\}\rangle = \widetilde{V}$.   \eop

\begin{lemma}\label{XXX.4}
Let $q$ be odd and greater than $3$. Then $K\subseteq \langle\tilde{\varepsilon}(l)\rangle$ for every hyperbolic line $l$. Moreover, $G$ acts trivially on $K$.
\end{lemma}
Let $l_0$ and $l_0'$ be two skew isotropic lines of $\Gamma$ and let $R$ be the regulus formed by the isotropic lines of $\Gamma$ that meet either of $l_0$ and $l'_0$. Let $R^*$ be the regulus opposite to $R$. Then $l_0, l'_0\in R^*$ and all lines of $R^*\setminus\{l_0,l'_0\}$ are hyperbolic (because $q$ is odd). Moreover $R^*\setminus\{l_0,l'_0\}$ can be partitioned in $(q-1)/2$ pairs of mutually orthogonal hyperbolic lines, say $\{l_1,l'_1\}$, $\{l_2,l'_2\},...,\{l_{(q-1)/2},l'_{(q-1)/2}\}$. We can choose distinct points $a_1,a_2,a_3\in l_0$, $a'_1, a'_2, a'_3\in l'_0$ and $b_1,..., b_d\in l_1$ such that $\langle \tilde{\varepsilon}(a_1), \tilde{\varepsilon}(a_2), \tilde{\varepsilon}(a_3)\rangle = \langle \tilde{\varepsilon}(l_0)\rangle$, $\langle \tilde{\varepsilon}(a'_1), \tilde{\varepsilon}(a'_2), \tilde{\varepsilon}(a'_3)\rangle = \langle \tilde{\varepsilon}(l'_0)\rangle$ and $\langle \tilde{\varepsilon}(b_1),..., \tilde{\varepsilon}(b_{d})\rangle = \langle \tilde{\varepsilon}(l_1)\rangle$. On the other hand, if $m\in R$ and $c_0 = m\cap l_0$, $c'_0 = m\cap l'_0$ and $c_1 = m\cap l_1$ then $\langle \tilde{\varepsilon}(c_0), \tilde{\varepsilon}(c'_0), \tilde{\varepsilon}(c_1)\rangle = \langle \tilde{\varepsilon}(m)\rangle$. It follows that $\langle \tilde{\varepsilon}(\{a_i\}_{i=1}^3\cup\{a'_i\}_{i=1}^3\cup\{b_i\}_{i=1}^d)\rangle = \langle \tilde{\varepsilon}(Q)\rangle$ where $Q := \cup_{m\in R}m = \cup_{l\in R^*}l$. Consequently, $\mathrm{dim}(\langle\tilde{\varepsilon}(Q)\rangle) \leq 6+d$.

On the other hand, $\mathrm{dim}(\langle \tilde{\varepsilon}(Q)\rangle)$ can also be computed as follows. Given  two hyperbolic lines $l, m\in R^*$ (possibly $m = l^\perp$) we can generate each of $\langle\tilde{\varepsilon}(l)\rangle$ and $\langle\tilde{\varepsilon}(m)\rangle$ with $d$ points chosen from $\tilde{\varepsilon}(l)$ and $\tilde{\varepsilon}(m)$ respectively. Moreover $3$ points chosen from the conic $\tilde{\varepsilon}(l_0)$ are enough to span the plane $\langle\tilde{\varepsilon}(l_0)\rangle$. Since every line of $R$ meets each of $l$, $m$ and $l_0$ in distinct points, the above $2d+3$ points are sufficient to span $\langle\tilde{\varepsilon}(Q)\rangle$. On the other hand, $\langle\varepsilon(l_0)\rangle\cap \langle \varepsilon(l\cup m)\rangle = 0$, as one can easily check. Hence $\mathrm{dim}(\langle\tilde{\varepsilon}(Q)\rangle) = 2d+3-\delta_{l,m}$, where $\delta_{l,m} := \mathrm{dim}(\langle\tilde{\varepsilon}(l)\rangle\cap\langle\tilde{\varepsilon}(m)\rangle)$. Therefore $2d+3-\delta_{l,m}\leq 6+d$, namely $d\leq 3+\delta_{l,m}$. However, $d = 3 + \delta_0$, where $\delta_0 := \mathrm{dim}(K_l)$ and $K_l = K\cap\langle\tilde{\varepsilon}(l)\rangle$. (Note that $\delta_0$ does not depend on the choice
of the hyperbolic line $l$, since $G$ is transitive on the set of hyperbolic lines and stabilizes $K$.) It follows that $\delta_0\leq \delta_{l,m}$. However, $\langle\varepsilon(l)\rangle\cap\langle\varepsilon(m)\rangle = 0$. Hence $\langle\tilde{\varepsilon}(l)\rangle\cap\langle\tilde{\varepsilon}(m)\rangle \subseteq K$, namely $\langle\tilde{\varepsilon}(l)\rangle\cap\langle\tilde{\varepsilon}(m)\rangle\subseteq K_l\cap K_m$. Therefore $\delta_0 \geq \delta_{l,m}$. By this inequality and the inequality $\delta_0 \leq \delta_{l,m}$ previously proved we get $\delta_0 = \delta_{l,m}$. Hence $K_l = K_m$.

So far, we have proved that $K_l = K_m$ for any two hyperbolic lines $l, m \in R^*$. Let $\equiv_K$ be the equivalence relation on the set of hyperbolic lines where $l\equiv_K m$ when $K_l = K_m$. By way of contradiction, suppose that the relation $\equiv_K$ is not the trivial relation, namely it admits at least two classes. Then the classes of $\equiv_K$ are imprimitivity classes for the group $G$. Each of them contains at least $q-1$ hyperbolic lines, contributed by a regulus like $R^*$, and it is partitioned in pairs of mutually orthogonal hyperbolic lines. Since $q > 3$ by assumption, each class of $\equiv_K$ contains at least two pairs of mutually orthogonal hyperbolic lines. It follows that the stabilizer in $G$ of a pair $\{l, l^\perp\}$ of mutually orthogonal hyperbolic lines is non-maximal. However, such a stabilizer is indeed a maximal subgroup of $G$ (compare Kleidman and Liebeck \cite{KL}, table 3.5.C, ${\cal C}_2$). Therefore $\equiv_K$ has only one class, namely $K_l = K_m $ for any two hyperbolic lines $l$ and $m$.

Put $K_0 := K_l$, where $l$ is any hyperbolic line. By the above, $K_0$ does not depend on the choice of the hyperbolic line $l$.
Let $l$ and $m$ be two intersecting hyperbolic lines. Then $\langle \tilde{\varepsilon}(l)\rangle/K_0= \langle \varepsilon(l)\rangle$; $\langle \tilde{\varepsilon}(m)\rangle/K_0= \langle \varepsilon(m)\rangle$ and $\mathrm{dim}(\langle \varepsilon(l)\rangle \cap \langle \varepsilon(m)\rangle)=1$. Hence $\delta +1=\mathrm{dim}(\langle \tilde{\varepsilon}(l)\rangle \cap \langle \tilde{\varepsilon}(m)\rangle)=\mathrm{dim}(K_0) +1=\delta_0+1$. So, $\delta = \delta_0$.

We shall now prove that $K_0 = K$. Let $a$ be a point and $l$ and $m$ two hyperbolic lines contained in $a^\perp$. Note that $l\cap m$ must be a point. By the proof of Lemma \ref{XXX.2}, the subspace $\langle \{\tilde{\varepsilon}(a)\}\cup\tilde{\varepsilon}(l)\cup\tilde{\varepsilon}(m)\rangle$ is either equal to $\langle\tilde{\varepsilon}(a^\perp)\rangle$ or it is a hyperplane of $\langle\tilde{\varepsilon}(a^\perp)\rangle$. If the latter case occurs, then $\langle\tilde{\varepsilon}(n)\rangle\not\subseteq\langle \{\tilde{\varepsilon}(a)\}\cup\tilde{\varepsilon}(l)\cup\tilde{\varepsilon}(m)\rangle$, for any hyperbolic line $n$ contained in $a^\perp$ but not passing through the point $l\cap m$. On the other hand $\langle \tilde{\varepsilon}(n)\rangle$ is generated by $K_n = K_0$ together with the $\tilde{\varepsilon}$-images of any three points $a_1, a_2, a_3\in n$. We may choose $a_1, a_2, a_3$ different from the intersection of $n$ with the line through $a$ and $l\cap m$.
Moreover, $K_0 \subseteq\langle \{\tilde{\varepsilon}(a)\}\cup\tilde{\varepsilon}(l)\cup\tilde{\varepsilon}(m)\rangle$. Hence $\langle\tilde{\varepsilon}(n)\rangle\subseteq\langle \{\tilde{\varepsilon}(a)\}\cup\tilde{\varepsilon}(l)\cup\tilde{\varepsilon}(m)\rangle$, contradicting the assumption that $\langle \{\tilde{\varepsilon}(a)\}\cup\tilde{\varepsilon}(l)\cup\tilde{\varepsilon}(m)\rangle$ is a hyperplane of $\langle\tilde{\varepsilon}(a^\perp)\rangle$. Hence
$\langle \{\tilde{\varepsilon}(a)\}\cup\tilde{\varepsilon}(l)\cup\tilde{\varepsilon}(m)\rangle
=\langle\tilde{\varepsilon}(a^\perp)\rangle$. As a consequence, $\mathrm{dim}(\langle\tilde{\varepsilon}(a^\perp)\rangle) = 2d - \delta = 6+\delta$. By the proof of Lemma \ref{XXX.3}, $\mathrm{dim}(\widetilde{V}) = 3d-2\delta+1 = 10+\delta$. If follows that $\mathrm{dim}(K) = \delta = \mathrm{dim}(K_0)$, namely $K = K_0$.

It remains to prove that $G$ acts trivially on $K$. Let $l$ and $m$ be two mutually orthogonal hyperbolic lines and $G_{l,m}$ the stabilizer of $l$ and $m$ in $G$. Then $G = A\times B$ where $A\cong B\cong \mathrm{SL}(2,q)$, $A$ acts faithfully on $l$ and trivially on $m$ while $B$ acts faithfully on $m$ and trivially on $l$. Accordingly, $A$ acts trivially on
$\langle\tilde{\varepsilon}(m)\rangle$ and $B$ acts trivially on $\langle\tilde{\varepsilon}(l)\rangle$. As $K\subseteq \langle\tilde{\varepsilon}(l)\rangle\cap\langle\tilde{\varepsilon}(m)\rangle$, both $A$ and $B$ act trivially on $K$, namely $G_{l,m}$ acts trivially on $K$. Thus, $G$ acts unfaithfully on $K$. However $G$ is simple. Hence it must act trivially on $K$.  \eop

\bigskip

In view of the final proposition of this section, we need a lemma on the action of the group $\mathrm{SL}(2,q)$ on the projective line $\mathrm{PG}(1,q)$. We denote the points of $\mathrm{PG}(1,q)$ by $a_\infty$ and $a_\lambda$ where $a_\infty = \langle (1,0)\rangle$ and $a_\lambda = \langle (\lambda,1)\rangle$ for $\lambda\in\mathbb{F}_q$. Let $\mathrm{SL}(2,q)_{\{\infty,0\}}$, $\mathrm{SL}(2,q)_{\{\infty,0\}}$ and $\mathrm{SL}(2,q)_{\{\infty, 0, 1\}}$ be the stabilizers in $\mathrm{SL}(2,q)$ of the ordered pair $(a_\infty, a_0)$, the unordered pair $\{a_\infty,a_0\}$ and the triple $\{a_\infty, a_0, a_1\}$ respectively. Assuming $q$ odd, let $\mathbb{F}_q^\Box$ be the subgroup of square elements of the multiplicative group $\mathbb{F}^*_q$ of $\mathbb{F}_q$ and put $\mathbb{F}^\diamondsuit_q = \mathbb{F}^*_q\setminus\mathbb{F}^\Box_q$. Accordingly, let $A^\Box := \{a_\lambda\}_{\lambda\in\mathbb{F}^\Box_q}$ and $A^\diamondsuit := \{a_\lambda\}_{\lambda\in\mathbb{F}^\diamondsuit_q}$.

\begin{lemma}\label{XXX.5}
Let $q$ be odd. \\
(1) The sets $A^\Box$ and $A^\diamondsuit$ are orbits of $\mathrm{SL}(2,q)_{\infty,0}$.\\
(2) If $-1\in \mathbb{F}^\diamondsuit_q$ then $\mathrm{SL}(2,q)_{\{\infty,0\}}$ contains an element switching $A^\Box$ with $A^\diamondsuit$. Consequently, $\mathrm{SL}(2,q)$ is transitive on the set of pairs $\{\{a, a'\}, a''\}$ with $a, a', a''$ distinct points of $\mathrm{PG}(1,q)$. \\
(3) If $-1\in \mathbb{F}^\Box_q$ then $\mathrm{SL}(2,q)_{\{\infty,0,1\}}$ contains an element mapping some element of $A^\Box$ onto an element of $A^\diamondsuit$.
\end{lemma}
\pr Most of the previous claims are well known. A few details remain which are perhaps not so well known, but their proofs are easy. We leave them to the reader.  \eop

\bigskip

The following proposition finishes the proof of Theorem \ref{main theorem univ1}

\begin{prop}\label{XXX.6}
Let $q$ be odd and $K \neq 0$. Then $q = 3$.
\end{prop}
\pr By way of contradiction, let $q > 3$. Let $l$ be a hyperbolic line and $G_l$ its stabilizer in $G$. Put $L := \tilde{\varepsilon}(l)$ and consider the action of $G_l \cong \mathrm{SL}(2,q)$ on the space $S = \langle L\rangle$. By Lemma \ref{XXX.4}, $K\subseteq S$ and $G_l$ acts trivially on $K$. Modulo factorizing by a subspace of $K$ of codimension 1, we may assume that $\mathrm{dim}(K) = 1$. Let $e_0, e_1,..., e_{q+1}$ be non-zero vectors of $S$ such that
$\langle e_0\rangle = K$ and $\langle e_1\rangle,..., \langle e_{q+1}\rangle$ are the points of $L$. Let $S_0 := \langle  e_1, e_2, e_3\rangle$ and, for every $i = 1,..., q+1$ let $p_i = \langle e_0, e_i\rangle\cap S_0$. Let $C := \{p_i\}_{i=1}^{q+1}$. As the projection of $L$ from $K$ is a conic, $C$ is a conic in the plane $S_0$. Without loss of generality, we may assume that the points of $C$ are represented by the following vectors, where $(s,t)\in\mathbb{F}^2_q\setminus\{(0,0)\}$:
\begin{equation}\label{XXXeq1}
u(s,t) := s^2e_1 + t^2e_2 + st(e_3-e_1-e_2).
\end{equation}
Clearly, $u(s,0) = s^2e_1$, $u(0,t) = t^2e_2$ and $u(s,s) = s^2e_3$. According to (\ref{XXXeq1}), the points of $\tilde{\varepsilon}(l)$ are represented by the following vectors:
\begin{equation}\label{XXXeq2}
v(s,t) := s^2e_1 + t^2e_2 + st(e_3-e_1-e_2) + f(s,t)e_0
\end{equation}
for a suitable function $f:\mathbb{F}^2_q\setminus\{(0,0)\}\rightarrow\mathbb{F}_q$ where $f(s,0) = f(0,t) = f(s,s) = 0$ and $f(rs,rt) = r^2f(s,t)$. In particular, $f(-s,-t) = f(s,t)$. Let now $G_{l,1,2}$ be the stabilizer of $\langle e_1\rangle$ and $\langle e_2\rangle$ in $G_l$ and let $g\in G_{l,1,2}$. Then
\begin{equation}\label{XXXeq3}
\left\{\begin{array}{rcl}
g(e_1) & = & \alpha e_1,\\
g(e_2) & = & \beta e_2,\\
g(e_3) & = &  a^2e_1 + b^2e_2 + ab(e_3-e_1-e_2) + f(a,b)e_0,
\end{array}\right.
\end{equation}
for suitable non-zero scalars $\alpha, \beta, a$ and $b$. Moreover $g(e_0) = e_0$. As $G_l\cong \mathrm{SL}(2,q)$, which is perfect, the elements of $G_l$ are represented by matrices of determinant 1 with respect to the basis $\{e_0, e_1, e_2, e_3\}$ of $S$. By this fact we immediately see that $\alpha\beta ab = 1.$ By linearity,
\[g(v(s,t)) =  s^2g(e_1) + t^2g(e_2) + st(g(e_3)-g(e_1)-g(e_2)) + f(s,t)g(e_0).\]
Substituting $g(e_1)$, $g(e_2)$ and $g(e_3)$ with their expressions as in (\ref{XXXeq3}) and by easy manipulations we obtain the following:
\begin{equation}\label{XXXeq5}
\begin{array}{ll}
g(v(s,t)) = & (s^2\alpha + st(a^2-\alpha))e_1 + (t^2\beta + st(b^2-\beta))e_2 + \\
 {} & + st\cdot ab(e_3-e_1-e_2) + (st\cdot f(a,b) + f(s,t))e_0.
\end{array}
\end{equation}
On the other hand,
\begin{equation}\label{XXXeq6}
g(v(s,t)) = u^2e_1 + v^2e_2 + uv(e_3-e_1-e_2) + f(u,v)e_0
\end{equation}
for suitable scalars $u, v$. By comparing (\ref{XXXeq5}) with (\ref{XXXeq6}) we see that
\begin{equation}\label{XXXeq7}
\begin{array}{l}
u^2  =  s^2\alpha + st(a^2-\alpha), ~~ v^2  =  t^2\beta + st(b^2-\beta), ~~ uv  =  st\cdot ab,\\
\mbox{and}~~  f(u,v) = st\cdot f(a,b) + f(s,t).
\end{array}
\end{equation}
Substituting the first three equalities of (\ref{XXXeq7}) in the identity $(uv)^2 = u^2v^2$, by standard manipulations we eventually obtain the following:
\begin{equation}\label{XXXeq10}
ab = 1, ~~ \alpha = a^2, ~~ \beta = b^2.
\end{equation}
We leave the details for the reader. Note that the computations to do here are just the same that we must do in order to compute the stabilizer in $X$ of two points $p_1$ and $p_2$ of a non-singular conic $C$ of $\PG(2,q)$, where $X$ is a subgroup of the stabilizer of $C$ in $\mathrm{SL}(3,q)$ and $X \cong \mathrm{SL}(2,q)$. (Recall that $G_l \cong \mathrm{SL}(2,q)$.)

Substituting $\alpha$ and $\beta$ with $a^2$ and $b^2$ in the first two equations of (\ref{XXXeq7}) and recalling that $uv = st\cdot ab$ by the third equation of (\ref{XXXeq7}) we get that either $u = sa$ and $v = tb$ or $u = -sa$ and $v = -tb$. Since $v(s,t) = v(-s,-t)$ we may assume without loss of generalities that $u = sa ~ \mbox{and} ~ v = tb.$ Substituting $u$ and $v$ with $sa$ and $tb$ in the last equation of (\ref{XXXeq7}) we obtain the following:

\begin{equation}\label{XXXeq11}
f(sa,tb) = f(s,t) + stf(a,b).
\end{equation}
The point $\langle g(e_3)\rangle = \langle v(a,b)\rangle$ only depends on the ration $a/b$. We have $a/b = a^2$ by (\ref{XXXeq10}). Hence $\F_q^\Box$ is the set of values taken by $a/b$ as $v(a,b)$ ranges in the $G_{l,1,2}$-orbit of $e_3$. If $s/t\in \mathbb{F}^\Box_q$ then we can switch $(s,t)$ with $(a,b)$ in (\ref{XXXeq11}), thus obtaining the following:
\begin{equation}\label{XXXeq12}
f(sa,tb) = f(a,b) + abf(s,t).
\end{equation}
Therefore, by comparing (\ref{XXXeq11}) with (\ref{XXXeq12}):

\begin{equation}\label{XXXeq13}
f(s,t)(1-ab) = f(a,b)(1-st)
\end{equation}
for any pair $(s,t)$ with $s/t\in \mathbb{F}^\Box_q$. Clearly, we can replace $(s,t)$ with $(s/t,1)$ in (\ref{XXXeq13}), obtaining the following:

\begin{equation}\label{XXXeq14}
f(\frac{s}{t},1)(1-ab) = f(a,b)(1-\frac{s}{t}\cdot t).
\end{equation}
On the other hand, if we divide both members of (\ref{XXXeq13}) by $t^2$ and recall that $f(s,t)/t^2 = f(s/t,1)$ then we also get the following:
\begin{equation}\label{XXXeq15}
f(\frac{s}{t},1)(1-ab) = f(a,b)(\frac{1}{t^2} - \frac{s}{t}).
\end{equation}
Since $t^2$ can be chosen to be different from 1 in (\ref{XXXeq15}), by comparing (\ref{XXXeq15}) with (\ref{XXXeq14}) we see that $f(a,b) = 0$. Hence $f(sa,ta) = f(s,t)$ by (\ref{XXXeq11}). Summarizing: if $a/b\in\mathbb{F}^\Box_q$ then $f(a,b) = 0$ and
\begin{equation}\label{XXXeq16}
f(sa,tb) = f(s,t)
\end{equation}
for any pair $(s,t)\in\mathbb{\F}^2_q\setminus\{(0,0)\}$. It now follows from (\ref{XXXeq2}) that the projective plane $\langle e_1, e_2, e_3\rangle$ contains at least $2+(q-1)/2$ points of $L$, namely $p_1 = \langle e_1\rangle$, $p_2 = \langle e_2\rangle$ and the $(q-1)/2$ points $\langle v(a,b)\rangle$ for $a/b\in \mathbb{F}_q^\Box$.

Choose a non-square $\eta\in\mathbb{F}^\diamondsuit_q$. For every $\mu\in\mathbb{F}^\diamondsuit_q$ there exists an element $\lambda\in \mathbb{F}^\Box_q$ such that $\mu = \eta\lambda$. By (\ref{XXXeq16}),
$f(\mu,1) = f(\eta,1)$. If $f(\eta,1) = 0$ then $f(s,t) = 0$ for every pair $(s,t)\in \mathbb{F}^2_q\setminus\{(0,0)\}$. In this case (\ref{XXXeq2}) shows that
$\langle e_1, e_2, e_3\rangle$ contains the whole of $L$. This contradicts the hypothesis that $S = \langle L\rangle$ is 4-dimensional. Therefore $f(\eta,1) \neq 0$ and $\langle e_1, e_2, e_3\rangle\cap\tilde{\varepsilon}(l) = \{p_1, p_2\}\cup \{\langle v(\lambda,1)\rangle\}_{\lambda\in \mathbb{F}_q^\Box}$.

Clearly, the action of $G_l$ on $L$ is isomorphic to the action of $\mathrm{SL}(2,q)$ on $\mathrm{PG}(1,q)$. We can choose an isomorphism $\pi$ between these two actions in such a way that $\pi(a_\infty) = p_1$, $\pi(a_0) = p_2$ and $\pi(a_1) = p_3$, where $a_\infty$, $a_0$ and $a_1$ are as stated in the paragraph before Lemma \ref{XXX.5}. By claim (1) of Lemma \ref{XXX.5}, $\pi(A^\Box) = \{p_\lambda\}_{\lambda\in\mathbb{F}^\Box_q}$ and $\pi(A^\diamondsuit) = \{p_\lambda\}_{\lambda\in\eta\mathbb{F}^\diamondsuit_q}$, where $p_\lambda := \langle v(\lambda,1)\rangle.$ Let $-1$ be a square. By claim (3) of Lemma \ref{XXX.5}, the projective plane $\langle e_1, e_2, e_3\rangle$ contains the point $p_\mu$ for some $\mu\in\mathbb{F}^\diamondsuit_q$, contrary to what we have proved in the previous paragraph. Therefore $-1$ is not a square. By (2) of Lemma \ref{XXX.5} and since $\{p_\lambda\}_{\lambda\in\mathbb{F}^\Box_q}$ is contained in the projective plane $\langle e_1, e_2, e_3\rangle$, the set $\{p_\lambda\}_{\lambda\in\eta\mathbb{F}^\diamondsuit_q}$ is contained in the plane $\langle e_1, e_2, v(-1,1)\rangle$. Moreover, $G_l$ is transitive on the set of pairs $\{\{p,p'\}, p''\}$ for $p, p', p''$ distinct points of $L$. Clearly, it is also transitive on the set of unordered triples of points of $L$.

We can now define the following point-block structure $\Lambda$ on $L$: the points of $\Lambda$ are the points of $L$ and the blocks are the intersection $X\cap L$, for $X$ a projective plane of $\mathrm{PG}(3,q)$. Since $G_l$ is transitive on the set of triples of points of $L$, the structure $\Lambda$ is a $3$-$(q+1,k,1)$ design with $k = 2+(q-1)/2 = (q+3)/2$. Let $N$ be the number of blocks of $\Lambda$. Then
\[N = \frac{{{q+1}\choose 3}}{{{(q+3)/2}\choose 3}} = \frac{8q}{q+3}.\]
Since $N$ is an integer, $q+3$ divides $8q$. Namely $q+3$ divides $24$. Recalling that $q$ is a prime power and $q > 3$ by assumption, the above divisibility conditions force $q = 9$ or $q = 5$. However, $-1$ is a square in either of $\mathbb{F}_5$ and $\mathbb{F}_9$, while we have previously proved that $-1$ cannot be a square. We have reached a final contradiction.  \eop

\subsection{Proof of Theorem \ref{main theorem univ2}} \label{sec4.2}

Let $n=2$. As remarked in the previous subsection, $\varepsilon^{\mathrm{spin}}_2$ embeds $\Delta_2$ as $W(3,\F)$ in $V(4,\F)$. As in the previous subsection, we put $\Gamma = W(3,\F)$, $\nu = \nu_4$ is the quadric veronesean map from $V(4,\F)$ to $V(10,\F)$, the embedding $\iota:\Gamma\rightarrow V(4,\F)$ is the inclusion of $W(3,\F)$ in $\PG(3,\F)$ and $\varepsilon := \nu\circ \iota$, but now $\F$ is a possibly infinite field of characteristic 2. Clearly $\varepsilon$, regarded as an embedding of $\Delta_2$, is isomorphic to $\varepsilon^{\mathrm{vs}}_2$.

Let $\tilde{\iota}:\Gamma\rightarrow\widetilde{V}$ be the hull of the embedding $\iota$. Then a suitable basis $B = \{e_i\}_{i=1}^4\cup\{e_j\}_{j\in J}$ can be chosen in $\widetilde{V}$ such that $\tilde{\iota}(\Gamma)$ is the generalized quadrangle associated with a quadratic form $\tilde{q}$ of $\widetilde{V}$ as follows, where the coordinates are taken with respect to $B$:
\[\tilde{q}(x_1,x_2,x_3,x_4,\{x_j\}_{j\in J}) = x_1x_3 + x_2x_4 + \sum_{j\in J}\lambda_jx^2_j\]
for suitable scalars $\lambda_j\in \F$ such that the form $\sum_{j\in J}\lambda_jx^2_j$ induced by $\tilde{q}$ on $\langle e_j\rangle_{j\in J}$ is totally anisotropic. Needless to say, $J\neq \emptyset$, $J\cap\{1,2,3,4\} = \emptyset$ and $|J|+4 = \mathrm{dim}(\widetilde{V})$. If $\F$ is perfect then $|J| = 1$, otherwise $|J| > 1$ (see e.g. De Bruyn and Pasini \cite{DBP}; we warn that if $\F$ is non-perfect the set $J$ can even be infinite). In any case, $\mathrm{dim}(\widetilde{V})\geq 5$.

It is convenient to give the set $I := \{1,2,3,4\}\cup J$ a total ordering, say $\leq$. We do not make any assumptions on $\leq$ but the following: the ordering $\leq$ induces on the set $\{1,2,3,4\}$ its natural ordering and $i \leq j$ for any $j\in J$ and $i = 1, 2, 3, 4$.

Let $\widetilde{W}$ be the subspace of $\widetilde{V}\otimes\widetilde{V}$ generated by the vectors $v\otimes v$ and $v\otimes w + w\otimes v$ for $v, w\in \widetilde{V}$. The set
\[B_{\mathrm{ver}} := \{e_i\otimes e_i\}_{i\in I}\cup\{e_i\otimes e_j + e_j\otimes e_i\}_{i,j\in I, i < j}\]
is a basis of $\widetilde{W}$. Let $\tilde{\nu}$ be the quadric veronesean map of $\widetilde{V}$ in $\widetilde{W}$, defined with respect to the bases $B$ and $B_{\mathrm{ver}}$ of $\widetilde{V}$ and $\widetilde{W}$. Then $\tilde{\varepsilon} := \tilde{\nu}\circ\tilde{\iota}$ is a veronesean embedding of $\Gamma$ in the hyperplane $H$ of $\widetilde{W}$ described by the following equation:
\[x_{1,3} + x_{2,4} + \sum_{i\in J}\lambda_jx_{j,j} = 0.\]
(Needless to say, coordinates are taken with respect to $B_{\mathrm{ver}}$.) It follows that $\mathrm{dim}(\tilde{\varepsilon}) = \mathrm{dim}(H) = \mathrm{dim}(\widetilde{W})-1$ ($= \mathrm{dim}(\widetilde{W})$ when the latter is infinite). The projection $\pi:\widetilde{V}\rightarrow V(4,\F)$ of $\tilde{\iota}$ onto $\iota$ naturally lifts to a linear mapping $\pi_{\mathrm{ver}}:\widetilde{W}\rightarrow V(10,\F)$ which maps $H$ onto $V(10,\F)$. Also, $\pi_{\mathrm{ver}}\circ\tilde{\varepsilon} = \varepsilon$. Hence $\pi_{\mathrm{ver}}$ is a morphism from $\tilde{\varepsilon}$ to $\varepsilon$. As $\mathrm{dim}(\tilde{\varepsilon}) = \mathrm{dim}(\widetilde{W})-1 \geq 14 > 10 = \mathrm{dim}(\varepsilon)$, the embedding $\varepsilon$ is not relatively universal.

\subsection{Proof of Theorem \ref{main theorem univ3}}

Let $n > 2$ and let $\F$ be a perfect field of characteristic 2. Then the spin embedding $\varepsilon^{\mathrm{spin}}_n$ is not relatively universal. Indeed, since $\F$ is a perfect field of characteristic 2 the building $\Delta$ can also be regarded as a building of type $C_n$. Accordingly, $\Delta_n$ also admits the embedding $\varepsilon^{\mathrm{sp}}_n$ defined at the end of Subsection \ref{1.2}. As remarked there, the embedding $\varepsilon_n^{\mathrm{spin}}$ is a proper quotient of $\varepsilon^{\mathrm{sp}}_n$ and $\mathrm{dim}(\varepsilon_n^{\mathrm{sp}}) = {{2n}\choose n}-{{2n}\choose{n-2}}$.

Put $\overline{V}_n := V({{2n}\choose n}-{{2n}\choose{n-2}},\F)$ and let $\overline{W}_n$ be the subspace of $\overline{V}_n\otimes\overline{V}_n$ generated by the vectors $v\otimes v$ and $v\otimes w + w\otimes v$ for $v, w\in \overline{V}_n$. Let $\bar{\nu}$ be the quadric veronesean map from $\overline{V}_n$ to $\overline{W}_n$ and put $\bar{\varepsilon} := \bar{\nu}\circ \varepsilon^{\mathrm{sp}}_n$. As in the previous subsection, one can show that $\varepsilon^{\mathrm{vs}}_n$ is a quotient of $\bar{\varepsilon}$.

Let now $X$ be an $(n-2)$-element of $\Delta$ and $Q_X$ the quad of $\Delta_n$ formed by the elements of $\Delta$ of type $n$ and $n-1$ incident to $X$. Then $\varepsilon_n^{\mathrm{sp}}(Q_X)$ is a generalized quadrangle of orthogonal type in a 5-subspace of $\overline{V}$. Therefore, as shown in the previous subsection,
 the veronesean embedding induced by $\bar{\varepsilon}$ on $Q_X$ is 14-dimensional. On the other hand, $\varepsilon_n^{\mathrm{spin}}(Q_X)$ is a generalized quadrangle of symplectic type in a 4-subspace of $V(2^n,\F)$. Consequently, $\varepsilon_n^{\mathrm{vs}}$ induces a 10-dimension embedding on $Q_X$. It follows that $\varepsilon_n^{\mathrm{vs}}\not\cong \bar{\varepsilon}$. Therefore $\varepsilon_n^{\mathrm{vs}} < \bar{\varepsilon}$. Theorem \ref{main theorem univ3} is proved.

\bigskip

\noindent
Ilaria Cardinali and Antonio Pasini,\\
Department of Information Engineering and Mathematics,\\
University of Siena,\\
Via Roma 56, I-53100, Siena, Italy\\
ilaria.cardinali@unisi.it, antonio.pasini@unisi.it
\end{document}